\documentclass[11pt,a4paper]{article} %  ÄÖÜ äöü ß
\usepackage[latin1]{inputenc}

\usepackage{amsmath}
\usepackage{amsfonts}
\usepackage{amssymb}
\usepackage{MnSymbol}
\usepackage{stmaryrd}

\newtheorem{satz}{Satz}[section]
\newtheorem{theorem}[satz]{Theorem}
\newtheorem{lemma}[satz]{Lemma}
\newtheorem{prop}[satz]{Proposition}
\newtheorem{cor}[satz]{Corollary}
\newtheorem{defin}[satz]{Definition}
\newtheorem{example}[satz]{Example}

\newcommand{\abs}[1]{\left|{#1}\right|}
\newcommand{\norm}[1]{\abs{\abs{{#1}}}}
\newcommand{\rund}[1]{\left(#1\right)}
\newcommand{\spitz}[1]{\left\langle{#1}\right\rangle}
\newcommand{\eckig}[1]{\left[{#1}\right]}
\newcommand{\schweif}[1]{\left\{#1\right\}}

\def\nz{{\rm I\kern-.20em N}}
\def\rz{{\rm I\kern-.20em R}}
\def\zz{\mathbb{Z}}

\def\H{\mathcal{H}}
\def\B{\mathcal{B}}
\def\C{\mathcal{C}}
\def\M{\mathcal{M}}
\def\N{\mathcal{N}}
\def\T{\mathcal{T}}
\def\I{\mathcal{I}}
\def\S{\mathcal{S}}
\def\P{\mathcal{P}}
\def\V{\mathcal{V}}
\def\X{\mathcal{X}}
\def\Y{\mathcal{Y}}
\def\Z{\mathcal{Z}}
\def\W{\mathcal{W}}
\def\Q{\mathcal{Q}}
\def\E{\mathcal{E}}
\def\F{\mathcal{F}}
\def\R{\mathcal{R}}
\def\L{\mathcal{L}}
\def\G{\mathcal{G}}

\def\bk{\mathbf{k}}
\def\br{\mathbf{r}}
\def\bs{\mathbf{s}}

\def\m{\mathfrak m}
\def\n{\mathfrak n}
\def\g{\mathfrak g}

\def\fX{\mathfrak X}

\def\Ker{{\rm Ker}}
\def\End{{\rm End}}
\def\triv{{\rm triv}}
\def\Hom{{\rm Hom}}
\def\id{{\rm id}}
\def\Id{{\rm Id}}
\def\der{{\rm der}}
\def\infaff{{\rm infaff}}
\def\Kill{{\rm Kill}}

\def\eps{\varepsilon}

\begin{document}

\vskip 1.0 true cm

\begin{center}

{\Huge \bf Integrating $\mathcal{P}$- super vectorfields}

\vskip 0.4 true cm

{\Huge \bf and the super geodesic flow}

\vskip 1.0 true cm

Roland Knevel \\

Bar-Ilan University RAMAT GAN

Department of Mathematics 

ISRAEL
\end{center}

\vskip 2.0 true cm



\begin{quote}
{\bf Abstract:} Aim of this article is to introduce the notion of integral and geodesic flows on $\mathcal{P}$-supermanifolds as certain partial actions of $\rz$~. First I introduce the concept of parametrization over a `small' super algebra $\mathcal{P}$~, which leads to the notion of $\mathcal{P}$-objects and is superized local deformation theory. It is shown how parametrization makes the theory much easier. A version of Palais' theorem for $\mathcal{P}$-supermanifolds is obtained stating that every infinitesimal $\mathcal{P}$-action of a simply connected $\mathcal{P}$- super Lie group $\G$ on a $\mathcal{P}$-supermanifold can be integrated to a whole action of $\G$~. Furthermore the faithful linearization of affine $\mathcal{P}$-supermorphisms is proven. Finally I show that Newton's, Lagrange's and Hamilton's approach to mechanics can be formulated also for $\mathcal{P}$- Riemannian supermanifolds and are infact equivalent.
\end{quote}

\section*{Introduction}

The motivation for this article comes from the general development of Riemannian super geometry, as introduced by O. Goertsches in his fundamental article \cite{Goertsches}, the applications to Riemannian symmetric superspaces and their importance in physics, see \cite{Goertsches} and re\-fe\-rences therein. Geodesics and integral flows are of course indispensable tools in differential geometry. A first question when trying to superize these notions to supermanifolds is: what is the appropriate counterpart to $\rz$ for supermanifolds? There is some belief that it is $\rz^{|1}$~, $\rz$ with one additional odd coordinate, and so super curves, supermorphisms from $I^{|1}$~, $I \subset \rz$ an open interval, into supermanifolds, and partial $\rz^{1|1}$-actions appear in the literature. However the benefit of these super curves and actions is questionable. First from a mathematical point of view: Have you ever observed what a nice and simple differential geometric object infact an open interval $I \subset \rz$ is?

\begin{itemize}
\item[(i)] All smooth vectorbundles on $I$ are trivial.
\item[(ii)] Every connection on a trivial vectorbundle on $I$ can be transformed into the trivial one by a smooth vectorbundle automorphism, and this is the reason for the existence of a parallel transport along every curve and of enough autoparallel curves, called geodesics.
\item[(iii)] Up to diffeomorphism there is only the trivial affine connection on $I$ - as can easily be shown with the help of the geodesic exponential map.
\item[(iv)] Given a $1$-dimensional connected Lie group $G$ with Lie algebra $\g$~, $t \mapsto \exp(t X)$ gives an isomorphism from either $(\rz, +)$ or $\rund{\rz / \zz, +}$ to $\G$~, $X$ an appropriately chosen generator of $\g$~.
\item[(v)] Every non-zero element of a Lie algebra generates a sub Lie algebra isomorphic to $\rz$~, and that is the reason for the existence of a Lie group exponential map and the possibility to integrate smooth vectorfields to partial actions of $\rz$~, called integral flows.
\end{itemize}

All these properties together explain the importance of curves and flows in super geometry. Unfortunately (ii), (iii), (iv) and (v) become wrong as soon as one adds odd coordinates to $\rz$ resp. $I$ : For (ii) and (iii) I will give counterexamples in this article myself; as a counterexample to (iv) in \cite{MontSan} J. Monterde and O. A. Sánchez-Valenzuela presented three non-isomorphic super Lie group structures on $\rz^{1|1}$ ; and that (v) does not hold in the context of super Lie groups is quite obvious. No wonder that severe problems occure when generalizing geodesics and integral flows to the super situation using the $\rz^{1|1}$-approach. For example in \cite{MontSan} the super integral flow, defined as a local super action of $\rz^{1|1}$~, does not exist to every super vectorfield $X$ due to the lack of a counterpart to (v): it exists iff the even and odd part of $X$ fulfill the same relations as the generators of the super Lie algebra of $\rz^{1|1}$. \\

Also from a physicists point of view the interpretation of super curves remains unclear: In most cases a curve in a manifold is interpreted as a physical state changing with time. But super symmetry does neither predict nor assume an additional odd time parameter. \\

On the other hand, given a supermanifold $\M$ with body $M$~, the ordinary curves in $\M$ cannot tell us anything about what is going on in the odd directions of $\M$ : since $\C^\infty_\rz$ has no non-zero nilpotent sections every curve $\gamma: I \rightarrow \M$ factors through its body map:

\[
\begin{array}{ccc}
M & \hookrightarrow & \M \\
\phantom{12}_{\gamma^\#} \nwarrow & \circlearrowleft & \nearrow {}_\gamma \phantom{12} \\
& I
\end{array}~.
\]

So using ordinary curves seems to be even worse. Nevertheless, super curves are still in general far away from separating the superfunctions on $\M$ : given a super curve $\gamma: I^{|1} \rightarrow \rz^{m|n}$~, $n \geq 2$~, and a non-zero function $f = f_I \xi^I$ on $\rz^{m|n}$ homogeneous in $\xi$ of degree $\geq 2$~, always $f \circ \gamma = 0$~. So what to do? \\

The solution seems to be parametrization, a superized local deformation theory: instead of single curves one should take whole families, parametrized by finitely many even and odd `parameters' generating a `small' super algebra. These parametrized curves turn out to separate the superfunctions, since already parametrized points do so! Geodesics in supermanifolds will be defined as such parametrized super curves in a quite natural way, and all problems disappear! Also super curves as described above obtain a nice interpretation in this concept: They are curves parametrized over $\bigwedge \rz$~, and the major difference to the above interpretation is that now one would never be inclined to study the derivative w.r.t. to the odd coordinate, which also from a physical point of view makes sense: super curves in physics should describe trajectories of super particles dealing with an even and odd state simultaneously. \\

But then the question is: why not parametrize everything? Supermanifolds, supermorphisms between them, super vectorbundles, connections, metrics, and super Lie group structures? We will of course do so! On the way we will see that supermanifolds and super vectorbundles are rigid, so admit {\bf no} non-trivial local super deformations and can up to isomorphism always be assumed to be unparametrized, while this is not true for the whole rest of the listed items! \\

Maybe highlights of this article are the faithful linearization result for affine $\P$-supermorphisms on connected $\P$-supermanifolds, corollary \ref{linearize}, and the generalization of Palais' theorem to the super case, theorem \ref{Palais}. The classical Palais' theorem, III of section IV.2 of \cite{Palais}, can be stated as follows:

\begin{theorem}[Palais' theorem]
Let $G$ be a simply connected Lie group with Lie algebra $\g$ (realized as the {\bf right}-invariant vectorfields on $G$ ), $M$ a smooth manifold and $\varphi: \g \rightarrow~\fX(M)$ a Lie algebra homomorphism, so an `infinitesimal action' of $G$~. Then $\varphi$ can be integrated to a whole smooth action of $G$ on $M$ iff $\varphi(\g)$ consists of complete vectorfields.
\end{theorem}

Throughout this article I use the ringed space description of supermanifolds as developped for example in \cite{Berezin} or \cite{DelMor} since it seems to be more adapted to the parametrization procedure. \\

{\bf Acknowledgement:} The research, which was finished - apart from the generalization of Palais' theorem - already in september 2010, was financed by the center of excellence grant of Israel Science Foundation (grant No.1691/10). I would like to thank the Bar-Ilan University, in particular Andre Reznikov, for the support of my research and scientific career, Oliver Goertsches from Hamburg for fruitfull discussions in november 2010 and Stephane Garnier from Metz for the suggestion to write super geodesics as integral curves to a super Hamiltonian vectorfield.

\section{$\mathcal{P}$-supermanifolds}

The body functor ${}^\#$ from the category of ringed spaces to the category of topological spaces associates to every ringed space $\X = (X, \S)$ its underlying topological space $\X^\# := X$ and to every morphism $\Phi: \X = (X, \S) \rightarrow \Y = (Y, \T)$ of ringed spaces its underlying continuous map $\Phi^\#: X \rightarrow Y$~. Given a ringed space $\X = (X, \S)$ and $U \subset X$ open, we can construct the ringed space $\X|_U := \rund{U, \S|_U}$~. We denote by ${\rm\bf SuperRingedSpac}$ the subcategory of all super ringed spaces, so ringed spaces $\X = (X, \S)$ where $\S$ is a sheaf of associative unital ( $\zz_2$-)graded rings, together with all supermorphisms, so morphisms $\Phi: \X = (X, \S) \rightarrow \Y = (Y, \T)$ whose associated sheaf morphism $\Phi^*: \T \rightarrow \Phi^\#_* \S$ is unital and respects the $\zz_2$-grading. Hereby $\Phi^\#_*$ denotes the push forward of sheaves under $\Phi^\#$~. \\

For the whole article let $\P$ and $\Q$ be finite dimensional real unital ($\zz_2$-)graded commutative algebras with unital algebra projections ${}^{\#'}: \P \rightarrow \rz$ and largest ideal $\Ker {}^{\#'} \lhd \P$ resp. $\Q$~, which is nilpotent. We call such algebras small, and the classical example of such a small algebra is $\bigwedge \rz^n$~. Furthermore let $\rho: \P \rightarrow \Q$ be an even unital algebra homomorphism and $\m$ the largest ideal of $\P$~. $\m$ is automatically graded, and its odd part equals that of $\P_1$~. \\

Let us already now fix some more notation:

\begin{itemize}
\item[(i)] Throughout this article I use Einstein notation. $\abs{R} \in \zz_2$ always denotes the parity of the homogeneous object $R$~, and given a collection $\rund{R_1, \dots, R_n}$ of homogeneous objects, $\abs{i} := \abs{R_i}$ for all $i = 1, \dots, n$~. For a finite set $I$~, $\abs{I}$ denotes its cardinality reduced to $\zz_2$~.
\item[(ii)] Given a graded commutative algebra $\S$~, graded $\S$-modules $\M$ and $\N$ and $a \in \M$~, $\M = \M_0 \oplus \M_1$ and $a = a_0 + a_1$ denote the splittings into even and odd components and $\M \boxtimes_\S \N$ the graded tensorproduct. We have of course a canonical isomorphism $\M \boxtimes \N \simeq \N \boxtimes \M \,~, \, a \otimes b \leftrightarrow (- 1)^{\abs{a} \abs{b}} b \otimes a$ for $a, b$ homogeneous, and for all $r, s \in \S$~, $a \in \M$ and $b \in \N$~, $a$ and $s$ homogeneous, we have $(r a) \otimes (s b) = (- 1)^{\abs{a} \abs{s}} r s (a \otimes b)$ in $\M \boxtimes \N$~.
\item[(iii)] As long as not stated the contrary, given a super ringed space $\X$~, $t$ denotes the canonical projection $\rz \times \X \twoheadrightarrow \rz$~.
\item[(iv)] The symbol $\diamondsuit$ means: evaluate this expression at the argument to obtain the desired map.

\end{itemize}

\begin{defin}
\item[(i)] A super ringed space $\M = \rund{M, \C^\infty_\M}$ such that $\C^\infty_\M$ is a sheaf of unital graded $\P$-algebras and $M = \M^\#$ is a smooth manifold of dimension $m$ is called a $\P$-supermanifold of super dimension $(m, n)$ iff locally $\C^\infty_\M \simeq \P \boxtimes \rund{\C^\infty_M \otimes \bigwedge \rz^n}$~.
\item[(ii)] Given $\P$-supermanifolds $\M$ and $\N$~, a supermorphism of ringed spaces $\Phi: \M \rightarrow \N$ is called a $\P$-supermorphism iff its associated sheaf morphism $\Phi^*: \C^\infty_\N \rightarrow \Phi^\#_* \C^\infty_\M$ is $\P$-linear.
\end{defin}

The $\P$-supermanifolds together with $\P$-supermorphisms form a subcategory $\P{\rm\bf -SuperMan}$ of ${\rm\bf SuperRingedSpac}$~. To $\rho: \P \rightarrow \Q$ we can associate a covariant functor ${}^\rho$ from $\P{\rm\bf -SuperMan}$ to $\Q{\rm\bf -SuperMan}$ assigning

\begin{itemize}
\item to a $\P$-supermanifold $\M = \rund{M, \C^\infty_\M}$ the $\Q$-supermanifold $\M^\rho := \rund{M, \Q \boxtimes_\P \C^\infty_\M}$ by using the multiplication $\Q \boxtimes \P \rightarrow \Q \,~, \, b \otimes a \mapsto b \rho(a)$~, and
\item to a $\P$-supermorphism $\Phi: \M \rightarrow \N$ the $\Q$-supermorphism $\Phi^\rho$ with $\rund{\Phi^\rho}^\# = \Phi^\#$ and

\[
\rund{\Phi^\rho}^*: \Q \boxtimes_\P \C^\infty_\N \rightarrow \Phi^\#_* \rund{\Q \boxtimes_\P \C^\infty_\M} = \Q \boxtimes_\P \Phi^\#_* \C^\infty_\M
\]

being the $\Q$ linear extension of $\Phi^*$~.
\end{itemize}

${}^\rho$ is clearly covariant functorial in $\rho$~, and ${}^\# \circ {}^\rho = {}^\#$~. In particular from the unital graded algebra homomorphisms

\[
\begin{array}{ccc}
\rz & \hookrightarrow & \P \\
\phantom{123}_{\id_\rz} \searrow & \circlearrowleft & \twoheadswarrow {}_{{}^{\#'}} \phantom{123} \\
& \rz &
\end{array}
\]

we obtain covariant functors

\[
\begin{array}{ccc}
{\rm\bf SuperMan} & \hookrightarrow & \P{\rm\bf -SuperMan} \\
\phantom{1234567890}_\Id \searrow & \circlearrowleft & \twoheadswarrow {}_{{}^{\#'}} \phantom{123456789012} \\
& {\rm\bf SuperMan} &
\end{array} \,~.
\]

The functor ${}^{\#'}$ is called the {\bf relative} body functor, {\bf not} to be mixed up with the body functor ${}^\#$ ! In the categrory of $\P$-supermanifolds we have a cross product: For $\P$-supermanifolds $\M$ and $\N$ with bodies $M$ and $N$ its cross product is given by $\M \times \N = \rund{M \times N, \rund{\Pr_M^{- 1} \C^\infty_\M} \boxtimes_\P \rund{\Pr_N^{- 1} \C^\infty_\N}}$~, where $\Pr_M^{- 1}$ and $\Pr_N^{- 1}$ denote the inverse image functors associated to the canonical projections $\Pr_M: M \times N \rightarrow M$ resp. $\Pr_N: M \times N \rightarrow N$~. \\

The local models of the $\P$-supermanifolds, which are used as local super charts, are the super open sets:

\[
U^{|n} := \rund{U, \C^\infty_{U^{|n}} } \,~,
\]

$U \subset \rz^m$ open and $\C^\infty_{U^{|n}} := \C^\infty_U \otimes \bigwedge \rz^n$~. They are ordinary supermanifolds but can be regarded as $\P$-supermanifolds as above by taking $\P \boxtimes \C^\infty_{U^{|n}}$ instead of $\C^\infty_{U^{|n}}$ as structure sheaf. On $U^{|n}$ we have the even coordinate functions $x^i \in \C^\infty(U) \hookrightarrow \C^\infty\rund{U^{|n}}_0$ and the odd ones $\xi^j \in \C^\infty\rund{U^{|n}}_1$~, which are just the standard base vectors of $\rz^n$~. Obviously, every section $f$ of $\P \boxtimes \C^\infty_{U^{|n}}$ can be decomposed as $f = f_I \xi^I$ with sections $f_I$ of $\P \otimes \C^\infty_U$~, $\xi^I = \xi^{j_1} \dots \xi^{j_q} \in \bigwedge \rz^n$~, $I = \schweif{j_1, \dots, j_q} \in \wp(\{1, \dots, n\})$~, $j_1 < \dots < j_q$~. Therefore we have the $\rz$-linear unital sheaf projection

\[
{}^\# := \P \boxtimes \C^\infty_{U^{|n}} \twoheadrightarrow \C^\infty_U \,~, \, f_I \xi^I \mapsto f_\emptyset^{\#'}
\]

and the $\C^\infty_{U^{|n}}$-linear extension 

\[
\begin{array}{ccc}
\P \boxtimes \C^\infty_{U^{|n}} & \mathop{\longrightarrow}\limits^{{}^{\rho}} & \Q \boxtimes \C^\infty_{U^{|n}} \\
\phantom{123}_{{}^\#} \twoheadsearrow & \circlearrowleft & \twoheadswarrow {}_{{}^\#} \phantom{12345} \\
& \C^\infty_U &
\end{array} \,~,
\]

of $\rho$~. By an easily done moderate generalization of Leites theorem, proposition 2.4 of \cite{DelMor}, we obtain

\begin{theorem} \label{Leites} For every $U \subset \rz^m$ and $V \subset \rz^p$ open there is a 1-1-corresponcence between $\P$-supermorphisms $\Phi: U^{|n} \rightarrow V^{|q}$ and tuples

\[
(f, \lambda) \in \rund{\P \boxtimes \C^\infty\rund{U^{|n}}}_0^{\oplus p} \oplus \rund{\P \boxtimes \C^\infty\rund{U^{|n}}}_1^{\oplus q}
\]

such that $f^\#(x) \in V$ for all $x \in U$~.

\begin{quote}
Obtaining the tuple $(f, \lambda)$ from $\Phi$ is easy:

\begin{equation} \label{Leitesform1}
(f, \lambda) = \Phi^* (y, \eta) \,~,
\end{equation}

where $y^k$ and $\lambda^l$ denote the super coordinate functions on $V^{|q}$~. \\

Obtaining $\Phi$ from the tuple $(f, \lambda)$ is more complicated: $\Phi$ is the unique $\P$-supermorphism such that $\Phi^\# = f^\#$ and for all $h = h_I \eta^I \in \P \boxtimes \C^\infty\rund{W^{|q}}$~, $W \subset V$ open, $h_I \in \P \otimes \C^\infty_V(W)$~,

\begin{equation} \label{Leitesform2}
\Phi^* h = \frac{1}{\br !} \rund{\rund{\partial_\br h_I} \circ f^\#} \rund{f - f^\#}^\br \lambda^I \in \C^\infty\rund{\rund{\rund{f^\#}^{- 1}}(W)}^{|n} \,~,
\end{equation}

where $\br$ runs through all multiindices in $\nz^m$ and is infact finite. This is of course a sort of Taylor formula.
\end{quote}
\end{theorem}

{\bf From now on for the rest of the article $\M$ and $\N$ denote $\P$-supermanifolds with bodies $M$ resp. $N$~.}

\begin{cor} \label{Leitescor}
\item[(i)] ${}^\#$ and ${}^\rho$ in local super charts glue together to unital $\rz$-linear even sheaf morphisms ${}^\#: \C^\infty_\M \twoheadrightarrow \C^\infty_M$~, whose kernel is the ideal of all nilpotent elements, and ${}^\rho: \C^\infty_\M \rightarrow \C^\infty_{\M^\rho}$~. ${}^\# \circ {}^\rho = {}^\#$~, and ${}^\rho$ is covariant in $\rho$~. $\rund{\Phi^* f}^\# = f^\# \circ \Phi^\#$ and $\rund{\Phi^* f}^\rho = \rund{\Phi^\rho}^* f^\rho$ for all $\P$-supermorphisms $\Phi: \M \rightarrow \N$~.

\item[(ii)] Given $V \subset \rz^p$ open, there exists a 1-1-correspondence between all $\P$-supermorphisms $\Phi: \M \rightarrow V^{|q}$ and tuples

\[
(f, \lambda) \in \C^\infty(\M)_0^{\oplus p} \oplus \C^\infty(\M)_0^{\oplus q}
\]

such that $f^\#(x) \in V$ for all $x \in M$~, in local super charts of $\M$ given by theorem \ref{Leites}. If $\Phi: \M \rightarrow V^{|q}$ is given by the tuple $(f, \lambda)$ and $\Xi: \N \rightarrow \M$ is a $\P$-supermorphism then $\Phi \circ \Xi$ is given by the tuple $\Xi^* (f, \lambda)$~, $\Phi^\#$ by $f^\#$ and $\Phi^\rho$ by $(f, \lambda)^\rho$~.
\end{cor}

{\it Proof:} (i) $\rund{\Phi^* f}^\# = f^\# \circ \Phi^\#$ and $\rund{\Phi^* f}^\rho = \rund{\Phi^\rho}^* f^\rho$ are correct as soon as $\Phi$ is a $\P$-supermorphism between super open sets, as it can be easily seen by formula (\ref{Leitesform2}).

(ii) obvious by (i) and formula (\ref{Leitesform1}). $\Box$ \\

Therefore from now on we identify a $\P$-supermorphism $\Phi: \M \rightarrow V^{|q}$ with its associated tuple $(f, \lambda) \in \C^\infty(\M)_0^{\oplus p} \oplus \C^\infty(\M)_0^{\oplus q}$ and write $f \circ \Phi := \Phi^* f$ heuristically thinking of `composing $\Phi$ with $f$ '.

Since a $\P$-supermanifold has no ordinary points, apart from the ones in its body, we have as a replacement the $\P$-points:

\begin{defin} A $\P$-supermorphism $x: \{0\} \rightarrow \M$ is called a $\P$-point of $\M$~. The set of all $\P$-points of $\M$ will be denoted by $\M^\P$~, and given a $\P$-supermorphism $\Phi: \M \rightarrow \N$ we write $\Phi^\P: \M^\P \rightarrow \N^\P \,~, \, x \mapsto \Phi(x) := \Phi \circ x$ thinking of `evaluating $\Phi$ at $x$ '.
\end{defin}

By theorem \ref{Leites} $\rund{U^{|n}}^\P$ is given by

\[
\rund{U^{|n}}^\P = \schweif{\left.(a, \beta) \in \P_0^{\oplus m} \oplus \P_1^{\oplus n} \, \right| \, a^\# \in U} = U \times \rund{\m_0^{\oplus m} \oplus \P_1^{\oplus n}} \subset \P_0^{\oplus m} \oplus \P_1^{\oplus n}
\]

open, and the body map ${}^\#: \rund{U^{|n}}^\P \twoheadrightarrow U$ is just the projection onto the first factor. Let $\Phi = (f, \lambda): U^{|n} \rightarrow V^{|q}$ be a $\P$-supermorphism, $V \subset \rz^p$ open. Then by corollary \ref{Leitescor} (ii) we see that $\Phi^\P: \rund{U^{|n}}^\P \rightarrow \rund{V^{|q}}^\P$ is smooth, given by

\[
(a, \beta) \mapsto \rund{\frac{1}{\bs !} \rund{\rund{\partial_\bs f_J}\rund{a^\#}} \rund{a - a^\#}^\br \beta^J \,~, \, \frac{1}{\bs !} \rund{\rund{\partial_\bs \lambda_J}\rund{a^\#}} \rund{a - a^\#}^\br \beta^J} \,~,
\]

and

\[
\begin{array}{ccc}
\phantom{12} \rund{U^{|n}}^\P & \mathop{\longrightarrow}\limits^{\Phi^\P} & \rund{V^{|q}}^\P \phantom{1} \\
{}^\# \downarrow & \circlearrowleft & \downarrow {}^\# \\
\phantom{1,} U & \mathop{\longrightarrow}\limits_{\Phi^\#} & V \phantom{1,}
\end{array} \,~.
\]

Therefore $\M^\P$ is not just a set, infact we obtain a covariant functor ${}^\P$ from $\P{\rm \bf -SuperMan}$ to the category of families $\Pi_X: X \twoheadrightarrow M$ of smooth manifolds together with smooth family morphisms

\[
\begin{array}{ccc}
\phantom{123,} X & \mathop{\longrightarrow}\limits^{\Psi} & Y \phantom{123} \\
\Pi_X \twoheaddownarrow & \circlearrowleft & \twoheaddownarrow \Pi_Y \\
\phantom{123,} M & \mathop{\longrightarrow}\limits_{\psi} & N \phantom{123,}
\end{array}
\]

assigning to every $\P$-supermanifold $\M$ its family of $\P$-points ${}^\#: \M^\P \twoheadrightarrow M$ and to every $\P$-supermorphism $\Phi: \M \rightarrow \N$ the smooth family morphism

\[
\begin{array}{ccc}
\phantom{12,} \M^\P & \mathop{\longrightarrow}\limits^{\Phi^\P} & \N^\P \phantom{1} \\
{}^\# \downarrow & \circlearrowleft & \downarrow {}^\# \\
\phantom{1,} M & \mathop{\longrightarrow}\limits_{\Phi^\#} & N \phantom{12}
\end{array}
\]

given by $\Phi^\P(x) := \Phi(x)$~. The obvious fundamental advantage of $\P$-points and difference to ordinary points is crucial:

\begin{lemma}
Let $\M$ be of super dimension $(m, n)$~. Then the $\rund{\P \boxtimes \bigwedge \rz^n}$-points of $\M$ separate the $\P$-supermorphisms $\M \rightarrow \N$~, more precisely: Let $\Phi, \Psi: \M \rightarrow \N$ be two $\P$-morphisms such that $\Phi(x) = \Psi(x)$ for all $x \in \M^{\P \boxtimes \bigwedge \rz^n}$~. Then $\Phi = \Psi$~.
\end{lemma}

Can I convince you that parametrization over $\P$ is superized local deformation theo\-ry? In the following description this becomes even more obvious: There is a contravariant functor $\widehat{\phantom{1}}$ from the category of small algebras together with unital even homomorphisms to ${\rm\bf SuperRingedSpac}$ assigning to every small algebra $\P$ the ringed space $(\{0\}, \P)$~. Obviously $\P{\rm\bf -SuperMan}$ is also the category of families of super ringed spaces $\Pi_\M: \M \twoheadrightarrow (\{0\}, \P)$ which locally look like

\[
\begin{array}{ccc}
\M & \mathop{\longrightarrow}\limits^{\sim} & \rz^{m|n} \times (\{0\}, \P) \\
\phantom{123}_{\Pi_\M} \twoheadsearrow & \circlearrowleft & \twoheadswarrow {}_{\Pr\nolimits_{(\{0\}, \P)} } \phantom{123456789} \\
& (\{0\}, \P) &
\end{array} \,~,
\]

and the functor ${}^\rho$ is nothing but the pullback under $\widehat \rho$~. In particular the canonical embedding ${\rm\bf SuperMan} \hookrightarrow \P{\rm\bf -SuperMan}$ means nothing but regarding a single object or morphism as a constant family over $(\{0\}, \P)$~, and ${}^{\#'}$ nothing but the restriction to the canonical embedding $\widehat{\phantom{1}^{\#'}}: \{0\} \hookrightarrow (\{0\}, \P)$~. The cross product of $\M$ and $\N$ just becomes the restricted cross product of families of ringed spaces $\M \times_{(\{0\}, \P)} \N$~, and finally the $\P$-points of $\M$ the $(\{0\}, \P)$-points of $\M$ from the functor of points approach to supermanifolds, see 2.8 and 2.9 of \cite{DelMor}. Since $\P$ is a local algebra we see that every $\P$-supermanifold $\M$ and every $\P$-supermorphism $\Phi$ between supermanifolds is a local deformation of its relative body $\M^{\#'}$ resp. $\Phi^{\#'}$~. On the other hand, as we know it from ordinary smooth manifolds, also a supermanifold does {\bf not} have non-trivial local deformations:

\pagebreak

\begin{theorem}[Rigidity of supermanifolds]
\item[(i)] There exists a $\P$-superdiffeomorphism $\Phi: \M^{\#'} \rightarrow \M$ with $\Phi^{\#'} = \Id_{\M^{\#'}}$~.
\item[(ii)] The group of $\P$-superdiffeomorphisms from $\M$ to itself acts transitively on the set of $\P$-points of $\M$~, in particular for every $x \in \M^\P$~, $\Phi$ in (i) can be chosen such that $\Phi\rund{x^\#} = x$~.
\end{theorem}

{\it Proof:} (i) As in the classical case by induction on the nilpontency degree of the maximal ideal $\m \lhd \P$ using the fact that the super tangent bundle $sT \M^{\#'}$ of $\M^{\#'}$ is a fine sheaf and so $H^1\rund{sT \M^{\#'}} = 0$~.

(ii) same as for smooth manifolds. $\Box$ \\

However, in general morphisms between supermanifolds {\bf do} have non-trivial local deformations:

\begin{example} Let $I \subset \rz$ be an open interval containing $0$~, $C \in \P_0 \setminus \{0\}$ such that $C^2 = 0$ and $\alpha \in \P_1 \setminus \{0\}$~. Let $\Phi$ be one of the following $\P$-supermorphisms

\[
C t : I \rightarrow I \phantom{1}~, \phantom{1} \alpha t: I \rightarrow \rz^{0|1} \phantom{1}~, \phantom{1} \alpha \xi: \rz^{0|1} \rightarrow I \phantom{1}~, \phantom{1} C \xi: \rz^{0|1} \rightarrow \rz^{0|1} \phantom{1}~.
\]

Then $\Phi^{\#'} = 0$~, and there exist {\bf no} $\P$-superdiffeomorphisms $\Omega$ and $\Xi$ of $I$ or $\rz^{0|1}$ appropriately such that $\Phi = \Xi \circ \Phi^{\#'} \circ \Omega$ since the right hand side would always be constant equal to $\Xi(0)$~.
\end{example}

Finally for later applications we need the following two results:

\begin{lemma} \label{homot}
$M$ and $\M^\P$ are homotopy equivalent.
\end{lemma}

{\it Proof:} Without restriction assume $\M$ is a supermanifold. There exists a canonical embedding $\iota: M \hookrightarrow \M^\P$ in local super charts of $\M$ given by $x \mapsto (x, 0)$~, and $\id_M = {}^\# \circ \iota$~. It remains to prove that $\id_{\M^\P}$ and $\iota \circ {}^\#$ are homotopic. For this purpose decompose $\m = \bigoplus_{r \in \nz \setminus \{0\}} \I^r$ such that all $\I^r$ are graded, $\m^r = \bigoplus_{k \geq r} \I^k$~, and so $\I^r \I^s \subset \bigoplus_{k \geq r + s} \I^k$ for all $r, s \in \nz \setminus \{0\}$~. Then we can decompose every $a \in \P$ as $a = a^{\#'} + \sum_{r \in \nz \setminus \{0\}} a_r$~, $a_r \in \I^r$~. So for every $c \in [0, 1]$ we obtain a unital even algebra homomorphism

\[
\rho_c: \P \rightarrow \P \,~, \, a = a^{\#'} + \sum_{r \in \nz \setminus \{0\}} a_r \mapsto a^{\#'} + \sum_{r \in \nz \setminus \{0\}} c^r a_r \,~,
\]

$c^r$ denoting the $r$-th power of $c$~. Obviously $\rho_1 = \id_\P$ and $\rho_0 = {}^{\#'}$~, which leads to a smooth map

\[
\sigma: [0, 1] \times \M^\P \rightarrow \M^\P \,~, \, x \mapsto x^{\rho_c}
\]

with $\sigma(1, \diamondsuit) = \id_{\M^\P}$ and $\sigma(0, \diamondsuit) = \iota \circ {}^\#$~. $\Box$

\begin{lemma}
Let $\Phi: \M \rightarrow \N$ be a $\P$-supermorphism. Then $\Phi$ is a $\P$-superdiffeomorphism iff $\Phi^{\#'}$ is a superdiffeomorphism.
\end{lemma}

{\it Proof:} `$\Rightarrow$': trivial.

`$\Leftarrow$': simple induction on the degree of nilpotency of the maximal ideal $\m \lhd \P$~. $\Box$ \\

So in particular the super inverse function theorem remains true also under parametrization.

\section{Basic differential geometry on $\mathcal{P}$- super vectorbundles} \label{vb}

Let $\X = (X, \S)$ be a super ringed space. Then we have the category ${\rm\bf grad}\S{\rm\bf -Mod}$ of graded $\S$-modules, and given a supermorphism $\Phi: \Y \rightarrow \X$~, $\Y = (Y, \T)$ a second super ringed space, the pullback $\Phi^*$ gives a covariant functor from ${\rm\bf grad}\S{\rm\bf -Mod}$ to ${\rm\bf grad}\T{\rm\bf -Mod}$~, contravariant in $\Phi$~. As for ordinary ringed spaces it is defined as the composition

\[
{\rm\bf grad}\S{\rm\bf -Mod} \mathop{\longrightarrow}\limits^{\rund{\Phi^\#}^{- 1}} {\rm\bf grad}\rund{\Phi^\#}^{- 1} \S{\rm\bf -Mod} \mathop{\longrightarrow}\limits^{\T \boxtimes_{\rund{\Phi^\#}^{- 1} \S}} {\rm\bf grad}\T{\rm\bf -Mod} \,~,
\]

where we use the multiplication

\[
\T \boxtimes \rund{\Phi^\#}^{- 1} \S \rightarrow \T \,~, \, f \otimes g \mapsto f \rund{\Phi^* g} \,~.
\]

Obviously $\Phi^*$ is compatible with taking graded tensor products, dual modules and transposition ${}^\dagger$~.

\begin{defin} Let $\E$ and $\F$ be graded $\S$-modules.
\item[(i)] $g \in (\E \boxtimes \E)^*_0$ is called a supermetric on $\E$ iff $g$ is graded symmetric, which means $\spitz{S \otimes T, g} = (- 1)^{\abs{S} \abs{T}} \spitz{T \otimes S, g}$ for all homogeneous sections $S$ and $T$ of $\E$~, and non-degenerate, which means that the induced homomorphism $\E \rightarrow \E^* \,~, \, T \mapsto \spitz{\diamondsuit \otimes T, g}$ is infact an isomorphism. We will use the notation $g(S, T) := \spitz{S \otimes T, g}$~. If $g$ is a supermetric then the pair $(\E, g)$ is called a super Riemannian $\S$-module.
\item[(ii)] Let $(\E, g)$ and $(\F, h)$ be super Riemannian $\S$-modules and $\varphi \in \Hom(\E, \F)_0$~. $\varphi$ is called isometric w.r.t. $g$ and $h$ iff $(\varphi \otimes \varphi)^\dagger h = g$~, in other words iff $h(\varphi(S), \varphi(T)) = g(S, T)$ for all sections $S$ and $T$ of $\E$~.
\end{defin}

Obviously all super Riemannian $\S$-modules with isometric even homomorphisms form a ca\-te\-go\-ry ${\rm\bf gradRiem}\S{\rm\bf -Mod}$~, and $\Phi^*$ becomes a covariant functor from ${\rm\bf gradRiem}\T{\rm\bf -Mod}$ to ${\rm\bf gradRiem}\S{\rm\bf -Mod}$~, contravariant in $\Phi$~. \\

Recall the definition of the graded $\S$-module $\S^{m|n}$ : As $\S$-module it is equal to $\S^{m + n}$~, but given a section $\rund{e^1, \dots, e^m, f^1, \dots, f^m}$ of $\S^{m|n}$~, it is homogeneous of parity $\eps$ iff all $e^i$ are homogeneous of partity $\eps$ and all $f^j$ of parity $\eps + 1$~.

\begin{defin} The subcategory ${\rm\bf SuperVB}(\X)$ of all graded $\S$-modules $\E$ for which there exists $(m, n) \in \nz^2$ such that locally $\E \simeq \S^{m|n}$ is called the category of super vectorbundles on $\X$~. We call $(m, n)$ the super rank of $\E$ and denote by ${\rm\bf RiemSuperVB}(\X)$ the subcategory of $({\rm Riem-}\S{\rm -gradMod})$ of all Riemannian super vectorbundles over $\X$~.
\end{defin}

Since $\Phi^* \S^{m|n} = \T^{m|n}$~, $\Phi^*$ maps super vectorbundles to super vectorbundles of same super rank. Let $\rund{e_k}$~, $\rund{f_l}$~, $\rund{g_k}$ and $\rund{h_l}$ denote the graded standard frames of $\S^{m|n}$~, $\S^{p|q}$~, $\T^{m|n}$ and $\T^{p|q}$ respectively. Then given a homomorphism $\varphi: \S^{m|n} \rightarrow \S^{p|q}$ by $\varphi\rund{e_k} = a_k^l f_l$~, $a_k^l \in \S(X)_{\abs{k} + \abs{l}}$~, its pullback $\Phi^* \varphi: \T^{m|n} \rightarrow \T^{p|q}$ is given by $\rund{\Phi^* \varphi} \rund{g_k} = \rund{\Phi^* a_k^l} f_l$~. Given a section $S = S^k e_k$~, $S^k \in \S(X)$~, of $\S^{m|n}$~, its pullback under $\Phi$ is given by $\Phi^* S = \rund{\Phi^* S^k} g_k$~. \\

Recall that, given super vectorbundles $\E$ and $\F$ on $\X$~, there exists a canonical isomorphism $\rund{\E \boxtimes \F}^* \simeq \E^* \boxtimes \F^*$ given by the pairing

\[
\spitz{S \otimes T, \alpha \otimes \beta} := (- 1)^{\abs{T} \abs{\alpha}} \spitz{S, \alpha} \spitz{T, \beta}
\]

for all sections $S, T, \alpha$ and $\beta$ of $\E$~, $\F$~, $\E^*$ and $\F^*$ resp., $T$ and $\beta$ homogeneous. \\

From now on we study $\P$- super vectorbundles on $\P$-supermanifolds.

\begin{defin} The category $\P{\rm\bf -superVB}(\M) := {\rm\bf superVB}\rund{M, \P \boxtimes \C^\infty_\M}$ is called the ca\-te\-go\-ry of $\P$- super vectorbundles on $\M$~.
\end{defin}

We have the body functor ${}^\#$ from $\P{\rm\bf -superVB}(\M)$ to the category ${\rm\bf VB}(M)$ of smooth vectorbundles over $M$ assigning to every $\P$-supervectorbundle $\E$ of super rank $(m, n)$ on $\M$ the vectorbundle $\E^\# := \rund{\C^\infty_M \otimes_{\C^\infty_\M} \E}_0$ of rank $m$ and to every homomorphism $\varphi: \E \rightarrow \F$ the restriction of its $\C^\infty_M$-linear extension to $\E^\#$~, which automatically maps to $\F^\#$ since $\C^\infty_M$ is purely even. It is covariant, compatible with pullbacks and commutes with taking the dual bundle and the graded tensor product.

We also have the covariant functor ${}^\rho := \C^\infty_{\M^\rho} \boxtimes_{\C^\infty_\M}$ from $\P{\rm\bf -superVB}(\M)$ to $\Q{\rm\bf -superVB}(\M)$~. It is covariant in $\rho$~, compatible with pullbacks, and ${}^\# \circ {}^\rho = {}^\#$~. Given a $\P$- super vectorbundle $\E$ on $\M$ of super rank $(p, q)$~, a section $S$ of $\E$ and $x \in \M^\P$~, we write $\E_x := x^* \E$ and think of it as the `fibre of $\E$ sitting at $x$ ', which is isomorphic to the graded $\P$-module $\P^{p|q}$~, and $S(x) \varphi_x := x^* S \in \E_x$~.

\begin{theorem}[rigidity of super vectorbundles] \label{noparametrizedVB}
Let $\E$ be a $\P$- super vectorbundle on $\M$~. Then there exists an isomorphism $\varphi: \E^{\#'} \mathop{\rightarrow}\limits^{\sim} \E$ with $\varphi^{\#'} = \id_{\E^{\#'}}$~.
\end{theorem}

{\it Proof:} As in the classical case by induction on the nilpotency degree of the maximal ideal $\m \lhd \P$ using that $\End \E^{\#'}$ is a fine sheaf and therefore $H^1\rund{\End \E^{\#'}} = 0$~. $\Box$ \\

Apart from the trivial bundles $\rund{\C^\infty_\M}^{r|s}$~, the most important $\P$-super vectorbundles on $\M$ are infact the super tangent bundle $sT \M$ and its dual bundle $sT^* \M$ : \\

$sT \M$ is the graded $\C^\infty_\M$-module of all $\P$-linear super derivations on $\C^\infty_\M$~. If $\M$ is of superdimension $(m, n)$~, $sT \M$ is of super rank $(m, n)$ since a local super chart on $\M$ gives a local frame $\rund{\partial_i}$ of $sT \M$ consisting of the partial derivatives w.r.t. the $m$ even and $n$ odd local coordinates $x^i$~. Infact the super commutator turns $sT \M$ into a sheaf of $\P$- super Lie algebras. The sections of $sT \M$ are called $\P$- super vectorfields, and as usual we write $\fX(\M)$ for the $\P$- super Lie algebra of global $\P$- super vectorfields.

$sT^* \M$ is also called the super cotangent bundle or the $\C^\infty_\M$-module of $1$-forms on $\M$~, and we have an even $\P$-linear sheaf morphism

\[
d: \C^\infty_\M \rightarrow sT^* \M
\]

given by $\spitz{X, d f} = X f$ for all $\P$- super vectorfields $X$ on $\M$ and smooth superfunctions $f$ on $\N$~. We have the Leibniz rule $d (f g) = (d f) g + f (d g)$ for all smooth superfunctions $f, g$ on $\M$~. In a local super chart of $\M$ obviously $\rund{d x^i}$ is the dual frame to $\rund{\partial_i}$~. \\

A $\P$-supermorphism $\Phi: \M \rightarrow \N$ induces an even $\P$-linear sheaf morphism $d \Phi$ from $sT \M$ into the graded $\C^\infty_\M$-module $s\der_\Phi\rund{\C^\infty_\N, \C^\infty_\M}$ of mixed $\P$- super derivations given by $((d \Phi) X) h := X (h \circ \Phi)$ for all $\P$- super vectorfields $X$ on $\M$ and smooth super functions $h$ on $\N$~. Fix a local super chart of $\N$~. Then taking the partial derivatives $\partial_j \in sT \N$~, $\rund{\partial_j \diamondsuit} \circ \Phi$ becomes a frame of $s\der_\Phi\rund{\C^\infty_\N, \C^\infty_\M}$~, we have a canonical $\C^\infty_\M$-linear identification

\[
s\der_\Phi\rund{\C^\infty_\M, \C^\infty_\N} \simeq \Phi^* sT \N \,~, \, \rund{\partial_j \diamondsuit} \circ \Phi \leftrightarrow \Phi^* \partial_j \,~,
\]

and $d \Phi$ is given by $X \mapsto \rund{X \Phi^k} \rund{\Phi^* \partial_k}$~.

Given a $\P$-supermorphism $\Phi : (\rz \times \M)|_U \rightarrow \N$~, $U \subset M$ open, we define \\
$\dot \Phi := (d \Phi) \partial_t \in \rund{\Phi^* sT \N}_0$~.

Obviously $(sT \M)^\# = T M$ and $(sT \M)^\rho = sT \M^\rho$~, and all above mentionned notions are compatible with ${}^\#$ and ${}^\rho$~.

\begin{defin} Let $\Phi: \M \rightarrow \N$ be a $\P$-supermorphism. Then $\P$- super vectorfields $X$ and $Y$ on $\M$ resp. $\N$ are called $\Phi$-related iff $(d \Phi) X = \Phi^* Y$~.
\end{defin}

An easy calculation shows that this notion is covariant in $\Phi$ and compatible with the super commutator, ${}^\#$ and ${}^\rho$~.

\begin{defin}
\item[(i)] Let $\E$ be a $\P$- super vectorbundle on $\M$~. An even $\P$-linear sheaf morphism \\
$\nabla: \E \rightarrow sT^* \M \boxtimes \E$ is called a $\P$-connection on $\E$ iff it fulfills the Leibniz rule

\[
\nabla(f S) = d f \otimes S + f \nabla S
\]

for all super functions $f$ on $\N$ and sections $S$ of $\E$~.
\item[(ii)] A $\P$-connection $\nabla$ on $sT \M$ is called an affine $\P$-connection on $\M$~.
\end{defin}

For $\P = \rz$ we obtain the definition of a usual connection on $\E$ resp. an affine connection on $\M$~, see for example \cite{DelMor}, section 3.6, or \cite{Goertsches}, section 4.2. Of course, given a $\P$-connection $\nabla$ on the $\P$ -super vectorbundle $\E$~, we can associate to it a connection $\nabla^\#$ on $\E^\#$ by the formula

\[
\nabla^\# S^\# = (\nabla S)^\#
\]

for all sections $S$ of $\E$ and a $\Q$-connection $\nabla^\rho$ on $\E^\rho$ by $\Q$-linear extension. Locally, using a local super chart of $\N$ and a local frame $\rund{e_j}$ of $\E$~, there is a 1-1 correspondence between $\P$-connections $\nabla$ on $\E$ and tuples of functions

\[
\Gamma_{i j}^k \in \rund{\C^\infty(\M) \boxtimes \P}_{\abs{i} + \abs{j} + \abs{k}} \,~,
\]

the Christoffel symbols associated to $\nabla$~, given by $\nabla_{\partial_i} e_j = \Gamma_{i j}^k e_k$~. $\nabla^\#$ is given by the Christoffel symbols $\rund{\Gamma_{i j}^k}^\#$~, $i, j, k$ even, and $\nabla^\rho$ by $\rund{\Gamma_{i j}^k}^\rho$~.

\begin{prop} Let $\Phi: \M \rightarrow \N$ be a $\P$-supermorphism and $\E$ a $\P$- super vectorbundle on $\N$ with $\P$-connection $\nabla$~. Then there exists a unique $\P$-connection $\Phi^* \nabla$ on $\Phi^* \E$ such that

\[
\rund{\Phi^* \nabla}_X \rund{\Phi^* S} = \spitz{(d \Phi) X, \Phi^* (\nabla S)}
\]

for all $\P$- super vectorfields $X$ on $\M$ and sections $S$ of $\E$~.
\end{prop}

{\it Proof:} Take local super charts of $\M$ and $\N$ and a local frame $\rund{e_k}$ of $\E$~. Let $\widehat \Gamma_{i k}^l$ denote the Christoffel symbols of $\Phi^* \nabla$ in the local super chart of $\M$ w.r.t. the local frame $\rund{\Phi^* e_k}$ of $\Phi^* \E$~. Now the condition is obviously fulfilled iff

\begin{equation} \label{connpullback}
\widehat \Gamma_{i k}^l = \rund{\partial_i \Phi^j} \rund{\Gamma_{j k}^l \circ \Phi} \,~. \, \Box
\end{equation}

\begin{defin}
$\Phi^* \nabla$ is called the pullback of $\nabla$ under $\Phi$~.
\end{defin}

$\Phi^* \nabla$ is contravariant in $\Phi$ and compatible with ${}^\rho$ and ${}^\#$~, which can be easily seen by formula (\ref{connpullback}).

\begin{lemma} \label{pullbackconnright}
Let $\Phi: \M \rightarrow \N$ be a $\P$-supermorphism and $X$ and $Y$ $\Phi$-related $\P$- super vectorfields on $\M$ resp. $\N$~. Let $\E$ be a $\P$- super vectorbundle on $\N$ with $\P$-connection $\nabla$ and $S$ a section of $\E$~. Then

\[
\rund{\Phi^* \nabla}_X \rund{\Phi^* S} = \Phi^* \rund{\nabla_Y S} \,~.
\]

\end{lemma}

{\it Proof:} $\rund{\Phi^* \nabla}_X \rund{\Phi^* S} = \spitz{(d \Phi) X, \Phi^* (\nabla S)} = \spitz{\Phi^* Y, \Phi^* (\nabla S)} = \Phi^* \rund{\nabla_Y S}$~. $\Box$

\begin{lemma} \label{constructconn} Let $\nabla^\E$~, $\nabla^\F$ be $\P$-connections on the $\P$- super vectorbundles $\E$ resp. $\F$ on $\M$~.
\item[(i)] There exists a unique $\P$-connection $\nabla^{\E^*}$ on $\E^*$ such that

\[
d \spitz{S, \psi} = \spitz{\nabla^\E S, \psi} + \spitz{S, \nabla^{\E^*} \psi}
\]

for all sections $S$ of $\E$ and $\psi$ of $\E^*$~.
\item[(ii)] There exists a unique $\P$-connection $\nabla^{\E \boxtimes \F}$ on $\E \boxtimes \F$ such that

\[
\nabla^{\E \boxtimes \F} (S \otimes T) = \rund{\nabla^\E S} \otimes T + S \otimes \rund{\nabla^\F T}
\]

for all sections $S$ of $\E$ and $T$ of $\F$~.
\item[(iii)] Passing from $\nabla^\E$ to $\nabla^{\E^*}$ and from $\nabla^\E$ and $\nabla^\F$ to $\nabla^{\E \boxtimes \F}$ are compatible with each other, with pullbacks under $\P$-supermorphisms, ${}^\#$ and ${}^\rho$~.
\end{lemma}

{\it Proof:}  Let $\Gamma_{i j}^k$ and ${\Gamma'}_{i r}^s$ denote the Christoffel symbols of $\nabla^\E$ resp. $\nabla^\F$ with respect to a local super chart of $\N$ and local graded frames $\rund{e_j}$ of $\E$ resp. $\rund{f_r}$ of $\F$~.

(i) Let $\widetilde \Gamma_{i r}^s$ denote the Christoffel symbols of $\nabla^{\E^*}$ w.r.t. the the dual frame $\rund{e_k^*}$ of $\E^*$~. Then the condition is fulfilled iff

\[
\widetilde \Gamma_{i r}^s = (- 1)^{1 + \abs{s} + \abs{r} \abs{s}} \Gamma_{i s}^r \,~.
\]

(ii) Let $\widehat \Gamma_{i (j, r)}^{(k, s)}$ denote the Christoffel symbols of $\nabla^{\E \boxtimes \F}$ w.r.t. the the local frame $\rund{e_j \otimes e_r}$ of $\E \boxtimes \F$~. Then the condition is fulfilled iff

\[
\widehat \Gamma_{i (j, r)}^{(k, s)} = \Gamma_{i j}^k \delta_r^s + (- 1)^{(\abs{r} + \abs{s}) \abs{j}} \delta_j^k {\Gamma'}_{i r}^s \,~.
\]

(iii) Let $S, T, \alpha$ and $\beta$ be sections in $\E$~, $\F$~, $\E^*$ and $\F^*$ resp., $T$ and $\alpha$ homogeneous. Then

\begin{eqnarray*}
&& \spitz{\nabla^{\E \boxtimes \F} (S \otimes T), \alpha \otimes \beta} + \spitz{S \otimes T, \nabla^{\E^* \boxtimes \F^*} (\alpha \otimes \beta)} \\
&& \phantom{1} = \spitz{\rund{\nabla^\E S} \otimes T + S \otimes \rund{\nabla^\F T}, \alpha \otimes \beta} + \spitz{S \otimes T, \rund{\nabla^{\E^*} \alpha} \otimes \beta + \alpha \otimes \rund{\nabla^{\F^*} \beta}} \\
&& \phantom{1} = (- 1)^{\abs{T} \abs{\alpha}} \rund{\rund{\spitz{\nabla^\E S, \alpha} + \spitz{S, \nabla^{\E^*} \alpha}} \spitz{T, \beta} \phantom{\frac{,}{,}} + \phantom{\frac{,}{,}} \spitz{S, \alpha} \rund{\spitz{\nabla^\F T, \beta} + \spitz{T, \nabla^{\F^*} \beta}} } \\
&& \phantom{1} = (- 1)^{\abs{T} \abs{\alpha}} \rund{\rund{d \spitz{S, \alpha}} \spitz{T, \beta} + \spitz{S, \alpha} \rund{d \spitz{T, \beta}} } \\
&& \phantom{1} = (- 1)^{\abs{T} \abs{\alpha}} d \rund{\spitz{S, \alpha} \spitz{T, \beta}} \\
&& \phantom{1} = d \spitz{S \otimes T, \alpha \otimes \beta} \,~,
\end{eqnarray*}

and so $\nabla^{\E^* \boxtimes \F^*}$ fulfills the required condition on $\nabla^{\rund{\E \boxtimes \F}^*}$ from (i). The rest is obvious by formula (\ref{connpullback}). $\Box$

\begin{defin} Let $\E, \F$ be $\P$- super vectorbundles on $\M$ with $\P$-connections $\nabla^\E$ resp. $\nabla^\F$~.
\item[(i)] A section $S$ of $\E$ is called parallel w.r.t. $\nabla$ iff $\nabla S = 0$~.
\item[(ii)] $\varphi \in \Hom(\E, \F) = \E^* \boxtimes \F$ is called parallel w.r.t. $\nabla^\E$ and $\nabla^\F$ iff it is parallel w.r.t. the $\P$-connection $\nabla^{\E^* \boxtimes \F}$ given by lemma \ref{constructconn}, or in other words iff $\nabla^\F_X \varphi(S) = \varphi\rund{\nabla^\E_X S}$ for all $\P$- super vectorfields $X$ on $\M$ and sections $S$ of $\E$~.
\end{defin}

Since $\nabla$ is even, all sections parallel w.r.t. $\nabla$ form a graded sub-$\P$-module of $\E$~. \\
 
Obviously the $\P$- super vectorbundles on $\M$ with $\P$-connections together with parallel homomorphisms form a category $\P{\rm\bf -SuperVBConn}(\M)$~, and given a $\P$-supermorphism $\Phi: \M \rightarrow \N$~, the pullback $\Phi^*$ gives a covariant functor from to $\P{\rm\bf -SuperVBConn}(\N)$ to $\P{\rm\bf -SuperVBConn}(\M)$~, contravariant in $\Phi$~. Moreover, we obtain covariant functors ${}^\#$ and ${}^\rho$ from $\P{\rm\bf -SuperVBConn}(\M)$ to the category of smooth vectorbundles on $M$ with connections resp. $\Q{\rm\bf -SuperVBConn}(\M)$ compatible with pullback. ${}^\# \circ {}^\rho = {}^\#$~, and ${}^\rho$ is covariant in $\rho$~. \\

Now we will generalize some of the nice geometric properties of $\rz$~, stated in the introduction, which will in particular lead to the notion of parallel transport:

\begin{theorem} \label{superVBonintervals} Let $I \subset \rz$ be an open interval.
\item[(i)] Every $\P$- super vectorbundle $\E$ on $I^{|n}$ of super rank $(p, q)$ is isomorphic to the trivial bundle $\rund{\C^\infty_{I^{|n}} }^{p|q}$~.
\item[(ii)] For every $\P$-connection $\nabla$ on $\rund{\C^\infty_I}^{p|q}$ there exists a $\P$-automorphism of $\rund{\C^\infty_I}^{p|q}$ parallel w.r.t. the trivial connection and $\nabla$~.
\item[(iii)] Given a $\P$-connection $\nabla$ on $\rund{\C^\infty_I}^{p|q}$ and $t_0 \in I^\P$ there exists a unique map

\[
\tau_{t_0}: \P^{p|q} \rightarrow \rund{\P \otimes \C^\infty_I}^{p|q}
\]

such that $|_{t_0} \circ \tau_{t_0} = \id_{\P^{\oplus p|q}}$ and $\nabla\rund{\tau_{t_0} v} = 0$ for all $v \in \P^{p|q}$~. $\tau_{t_0}$ $\P$-linear and even.
\end{theorem}

{\it Proof:} (i) Obviously we have an equality of categories

\[
\P{\rm\bf -superVB}\rund{I^{|n}} = \rund{\P \otimes \bigwedge \rz^n}{\rm\bf -superVB}(I) \,~,
\]

so by theorem \ref{noparametrizedVB} we may assume without restriction that $\E$ is a super vectorbundle on $I$~. Then since even and odd part of $\E$ are ordinary vectorbundles on $I$~, the result follows from classical analysis.

(ii) Choose $t_0 \in I$~. Then given two $\P$-connections $\nabla$ and $\widetilde \nabla$ on $\rund{\C^\infty_I}^{p|q}$~, there exists a unique even homomorphism $\varphi^{\nabla, \widetilde \nabla}: \rund{\C^\infty_I}^{p|q} \rightarrow \rund{\C^\infty_I}^{p|q}$ parallel w.r.t. $\nabla$ and $\widetilde \nabla$ such that $\varphi^{\nabla, \widetilde \nabla}\rund{t_0} = \id_{\P^{p|q}}$~. Indeed: Let $\Gamma_k^l, {\widetilde \Gamma}_k^l \in \C^\infty(I) \otimes \P_{\abs{k} + \abs{l}}$ $\widetilde \Gamma_k^l$ denote the Christoffel symbols of $\nabla$ resp. $\widetilde \nabla$~, and write $\varphi^{\nabla, \widetilde \nabla}\rund{e_k} = a_k^l e_l$~, $a \in \rund{\C^\infty(I) \otimes \P}^{(m|n) \times (m|n)}_0$~. Then the condition on $\varphi^{\nabla, \widetilde \nabla}$ is equivalent to the initial value problem $a_k^r\rund{t_0} = \delta_k^r$ and

\[
\dot{a_k^r} = \Gamma_k^l a_l^r - a_k^l {\widetilde \Gamma}_l^r \,~,
\]

which by the Picard-Lindelöf theorem possesses exactly one solution $a$~.

Now by uniqueness we obtain $\varphi^{\nabla, \widetilde \nabla} \circ \varphi^{\nabla, \widetilde \nabla} = \id_\E$~, and so all $\varphi^{\nabla, \widetilde \nabla}$ are infact isomorphisms.

(iii) By (ii) we may assume that $\nabla$ is the trivial connection, and then the canonical embedding as constant sections is obviously the only possible choice for $\tau_{t_0}$~. $\Box$

\begin{defin}
\item[(i)] $\tau_{t_0}: \P^{\oplus p|q} \rightarrow \rund{\P \otimes \C^\infty_I}^{p|q}$ is called the parallel transport along $I$~.
\item[(ii)] Let $\gamma: I \rightarrow \M$ be a $\P$-curve, so a $\P$-supermorphism $\gamma: I \rightarrow \M$~, $I \subset \rz$ an open interval, and $\E$ a $\P$- super vectorbundle on $\M$ with $\P$-connection $\nabla$~. Then the even $\P$-linear map

\[
\nabla_{\dot \gamma} := \rund{\gamma^* \nabla}_{\partial_t} : \gamma^* \E \rightarrow \gamma^* \E
\]

is called the covariant derivative along $\gamma$~.
\end{defin}

Of course both notions are compatible with ${}^\#$ and ${}^\rho$~. Let $\E$ be a $\P$- super vectorbundle on $\M$~. Then after choosing a local frame $\rund{e_k}$ of $\E$ we can write any section $S$ of $\gamma^* \E$ uniquely as $S = S^k \rund{\gamma^* e_k}$~, $S^k$ smooth functions on $I$~. In a local super chart of $\M$ using the Christoffel symbols $\Gamma_{i k}^l$ of $\nabla$ and formula (\ref{connpullback}) we obtain

\[
\nabla_{\dot \gamma} S = \dot S + S^j \dot{\gamma^i} \rund{\Gamma_{i j}^k \circ \gamma} e_k \,~,
\]

which coincides with the definition of the covariant derivative along $\gamma$ in \cite{Goertsches}, 4.3, for $\P = \rz$~.

\begin{theorem}[Invariance of the covariant derivative under $\P$-supermorphisms] \label{covderpullback} Let $\gamma: I \rightarrow \M$ be a $\P$-curve and $\Phi: \M \rightarrow \N$ a $\P$-supermorphism.
\item[(i)] Let $\E$ be a $\P$- super vectorbundle on $\N$ with $\P$-connection $\nabla$~. Then $\nabla_{(\Phi \circ \gamma)^\cdot} = \rund{\Phi^* \nabla}_{\dot \gamma}$~.
\item[(ii)] Let $\Phi$ be affine w.r.t. $\nabla^\M$ and $\nabla^\N$~. Then

\[
\begin{array}{ccc}
\phantom{\gamma^* (d \Phi)} \gamma^* sT \M & \mathop{\longrightarrow}\limits^{\nabla_{\dot \gamma}^\M} & \gamma^* sT \M \phantom{\gamma^* (d \Phi)} \\
\gamma^* (d \Phi) \downarrow & \circlearrowleft & \downarrow \gamma^* (d \Phi) \\
\phantom{\gamma^* (d \Phi)} \gamma^* \Phi^* sT \N & \mathop{\longrightarrow}\limits_{\nabla_{(\Phi \circ \gamma)^\cdot}^\N} & \gamma^* \Phi^* sT \N \phantom{\gamma^* (d \Phi)}
\end{array} \,~.
\]

\end{theorem}

{\it Proof:} (i) just the contravariant functoriality of $\Phi^* \nabla$ w.r.t. $\Phi$~.

(ii) Let $S \in \gamma^* sT \M$~. Then since $d \Phi$ is parallel w.r.t. $\nabla^\M$ and $\Phi^* \nabla^\N$~, also $\gamma^* (d \Phi)$ is parallel w.r.t. $\gamma^* \nabla^\M$ and $\gamma^* \Phi^* \nabla^\N$~. Therefore

\begin{eqnarray*}
\rund{\nabla^\N_{(\Phi \circ \gamma)^\cdot} \circ \gamma^* (d \Phi)} S &=& \rund{\gamma^* \Phi^* \nabla^\N}_{\partial_t} \rund{\rund{\gamma^* (d \Phi)} S} \\
&=& \rund{\gamma^* (d \Phi)} \rund{\rund{\gamma^* \nabla^\M}_{\partial_t} S} \\
&=& \rund{\rund{\gamma^* (d \Phi)} \circ \nabla^\M_{\dot \gamma}} S \,~. \, \Box
\end{eqnarray*}

(ii) and (iii) of theorem \ref{superVBonintervals} are indeed false on $I^{|n}$~, even among unparametrized connections:

\begin{example}
Let $\E$ denote the trivial super vectorbundle on $\rz^{0|1}$ of super rank $(1, 0)$ or $(0, 1)$~, in other words $\E = \rund{\bigwedge \rz}^{1|0}$ resp. $\E = \rund{\bigwedge \rz}^{0|1}$ as graded $\bigwedge \rz$-module. Let $\xi$ denote the odd coordinate on $\rz^{0|1}$ and $e$ the standard frame on $\E$~. Define the family $\rund{\nabla^C}_{C \in \rz}$ of connections on $\E$ by $\nabla^C_{\partial_\xi} e = C \xi e$~. Since it is smooth we may also take $C \in \P_0$ and so obtain a $\P$-connection on $\E$~.

\begin{itemize}
\item[(i)] Let $C \in \P_0 \setminus \{0\}$~. Then there exists {\bf no} $\P$-section $S$ of $\E$ parallel w.r.t. $\nabla^C$ with $S(0) = e$~. Indeed: Write $S = (1 + b \xi) e$~, $b \in \P$~. Then $\nabla_{\partial_\xi} S = (b + C \xi) e \not= 0$~. So there there exists {\bf no} $\P$-isomorphism $\varphi: \E \rightarrow \E$ parallel w.r.t. the trivial connection and $\nabla$~. In particular for $C^2 = 0$ we obtain a non-trivial infinitesimal deformation of the trivial connection on $\E$~.
\item[(ii)] $\nabla^C$~, $C \in \P_0 \setminus \{0\}$~, is {\bf not} flat: Following \cite{Goertsches}, 4.2, we define for every $\P$-connection $\nabla$ of a $\P$- super vectorbundle its curvature tensor $R_\nabla \in \rund{sT^* \M \boxtimes sT^* \M \boxtimes \End \E}_0(M)$ as

\[
R_\nabla(X, Y) S := \eckig{\nabla_X, \nabla_Y} S - \nabla_{\eckig{X, Y}} S
\]

for all $\P$- super vectorfields $X, Y$ and sections $S$ of $\E$~. Here $R_{\nabla^C}$ is given by $R_{\nabla^C}\rund{\partial_\xi, \partial_\xi} = 2 C$~.
\item[(iii)] If $sT \ \rz^{0|1} \simeq \E$ via $\partial_\xi \leftrightarrow e$ then $\nabla^C$~, $C \in \P_0 \setminus \{0\}$~, is also {\bf not} torsionfree: Again following \cite{Goertsches}, 4.2, we define the torsion tensor $T_\nabla \in \rund{sT \M \boxtimes sT \M}^*_0(M)$ of an affine $\P$-connection $\nabla$ on $\M$ as

\[
T_\nabla(X, Y) := \spitz{T_\nabla, X \otimes Y} := \nabla_X Y - (- 1)^{\abs{X} \abs{Y}} \nabla_Y X - \eckig{X, Y}
\]

for all homogeneous $\P$- super vectorfields $X, Y$~. Here $T_{\nabla^C}\rund{\partial_\xi, \partial_\xi} = 2 C \xi \partial_\xi$~.
\end{itemize}
\end{example}

\begin{defin} Let $\Phi: \M \rightarrow \N$ be a $\P$-supermorphism.
\item[(i)] A supermetric on $sT \M$ is called a Riemannian $\P$-supermetric on $\M$~. The pair $(\M, g)$~, $g$ a  Riemannian $\P$-supermetric, will be called a $\P$- Riemannian supermanifold.
\item[(ii)] $\Phi$ is called isometric w.r.t. Riemannian $\P$-supermetrics $g$ on $\M$ and $h$ on $\N$ iff $d \Phi$ is isometric w.r.t. $g$ and $\Phi^* h$~.
\item[(iii)] $\Phi$ is called affine w.r.t. the affine $\P$-connections $\nabla^\M$ on $\M$ and $\nabla^\N$ on $\N$ iff $d \Phi$ is parallel w.r.t. $\nabla^\M$ and $\Phi^* \nabla^\N$~, in other words iff

\[
(d \Phi) \rund{\nabla^\M_X Y} = \rund{\Phi^* \nabla^\N}_X ((d \Phi) Y)
\]

for all $\P$- super vectorfields $X$ and $Y$ on $\M$~.
\end{defin}

\begin{theorem}
\item[(i)] The $\P$- Riemannian supermanifolds together with isometric $\P$-supermorphisms form a category $\P{\rm\bf -RiemSuperMan}$~.
\item[(ii)] The $\P$-supermanifolds with affine $\P$-connections together with affine $\P$-supermorphisms form a category $\P{\rm\bf -SuperManAffConn}$~.
\item[(iii)] ${}^\#$ and ${}^\rho$ give covariant functors from $\P{\rm\bf -RiemSuperMan}$ to the category of smooth Riemannian manifolds resp. $\Q{\rm\bf -RiemSuperMan}$ and from $\P{\rm\bf -SuperManAffConn}$ to the category of smooth manifolds with affine connections together with parallel smooth maps resp. $\Q{\rm\bf -SuperManAffConn}$~. ${}^\# \circ {}^\rho = {}^\#$~, and ${}^\rho$ is covariant in $\rho$~.
\end{theorem}

{\it Proof:} Let $\Phi: \M \rightarrow \N$ and $\Psi: \N \rightarrow \R$ be $\P$-supermorphisms.

(i) Let $\Phi$ and $\Psi$ be isometric w.r.t. the Riemannian $\P$-supermetrics $g^\M$~, $g^\N$ and $g^\R$~, which means $d \Phi$ and $d \Psi$ are isometric w.r.t. $g^\M$ and $\Phi^* g^\N$ resp. $g^\N$ and $\Psi^* g^\R$~. Then also $\Phi^* (d \Psi)$ is isometric w.r.t. $\Phi^* g^\N$ and $\Phi^* \Psi^* g^\R = (\Psi \circ \Phi)^* g^\R$~, and so $d (\Psi \circ \Phi) = \rund{\Phi^* (d \Psi)} d \Phi$ is isometric w.r.t. $g^\M$ and $(\Psi \circ \Phi)^* g^\R$~, which says that $\Psi \circ \Phi$ is isometric.

(ii) similar to (i).

(iii) obvious. $\Box$

\begin{lemma} Let $\Phi: \M \rightarrow \N$ be a $\P$-supermorphism and $X, Y$ be $\P$- super vectorfields on $\M$~, $Z, W$ on $\N$ such that $X$ and $Z$ and also $Y$ and $W$ are $\Phi$-related.
\item[(i)] If $\Phi$ is isometric w.r.t. the Riemannian $\P$-supermetrics $g^\M$ and $g^\N$ then \\
$g^\N(Z, W) \circ \Phi = g^\M(X, Y)$~.
\item[(ii)] If $\Phi$ is affine w.r.t. the affine $\P$-connections $\nabla^\M$ and $\nabla^\N$ then also $\nabla^\M_X Y$ and $\nabla^\N_Z W$ are $\Phi$-related.
\end{lemma}

{\it Proof:} (i) $g^\N(Z, W) \circ \Phi = \rund{\Phi^* g^\N}\rund{\Phi^* Z, \Phi^* W} = \rund{\Phi^* g^\N}\rund{(d \Phi) X, (d \Phi) Y} = g^\M(X, Y)$~.

(ii) We have to show that $(d \Phi)\rund{\nabla^\M_X Z} = \Phi^* \rund{\nabla^\N_Y W}$~. By lemma \ref{pullbackconnright}

\[
(d \Phi)\rund{\nabla^\M_X Z} = \rund{\Phi^* \nabla^\N}_X \rund{(d \Phi) Z} = \rund{\Phi^* \nabla^\N}_X \rund{\Phi^* W} = \Phi^* \rund{\nabla^\N_Y W} \,~. \, \Box
\]

As in \cite{Goertsches}:

\begin{defin} Let $(\E, g)$ be a Riemannian $\P$- super vectorbundle on $\M$ with $\P$-connection $\nabla$~. Then $\nabla$ and $g$ are called compatible iff $\nabla^{(\E \boxtimes \E)^*} g = 0$~, $\nabla^{(\E \boxtimes \E)^*}$ given by lemma \ref{constructconn} (i) and (ii), or in other words iff $X g(S, T) = g\rund{\nabla_X S, T} + (- 1)^{\abs{X} \abs{S}} g\rund{S, \nabla_X T}$ for all homogeneous vectorfiels $X$ on $\M$ and sections $S$~, $T$ of $\E$~, $S$ homogeneous.
\end{defin}

\begin{theorem}[Levi-Civita connection] \label{LeviCivita}
Let $g$ be a Riemannian $\P$-supermetric on $\M$~. Then there exists a unique torsionfree affine $\P$-connection $\nabla$~, called the Levi-Civita connection, on $\M$ compatible with $g$~. In a local super chart its Christoffel symbols are given by

\begin{equation} \label{LeviCivitaform}
2 \Gamma_{i j}^k g_{k r} = \partial_i g_{j r} + (- 1)^{\abs{i} \abs{j}} \partial_j g_{i r} - (- 1)^{\abs{r} (\abs{i} + \abs{j}) } \partial_r g_{i j}
\end{equation}

with $g_{i j} = g\rund{\partial_i, \partial_j} \in \C^\infty(\M)_{\abs{i} + \abs{j}}$~, $g_{j i} = (- 1)^{\abs{i} \abs{j}} g_{i j}$~.
\end{theorem}

By formula (\ref{LeviCivitaform}) it is obvious that passing from $g$ to its associated Levi-Civita connection commutes with ${}^\#$ and ${}^\rho$~. Observe that for an affine $\P$-connection $\nabla$ the following are equivalent:

\begin{itemize}
\item[(i)] $\nabla$ is torsionfree,
\item[(ii)] for all $\P$-supermorphisms $\Phi: U \rightarrow \M$~, $U \subset \rz^2$ open,

\[
\rund{\varphi^* \nabla}_{\partial_s} \rund{(d \Phi) \partial_t} = \rund{\varphi^* \nabla}_{\partial_t} \rund{(d \Phi) \partial_s} \,~,
\]

\item[(iii)] in local super charts the Christoffel symbols fulfill $\Gamma_{i j}^k = (- 1)^{\abs{i} \abs{j}} \Gamma_{j i}^k$~.
\end{itemize}

For later purpose we need:

\begin{lemma} \label{isometric} Let $\Phi: \M \rightarrow \N$ be an isometric $\P$-supermorphism between the $\P$- Riemannian supermanifolds $(\M, g)$ and $(\N, h)$~.
\item[(i)] $\Phi$ is an immersion.
\item[(ii)] If $\M$ and $\N$ are of the same super dimension then $\Phi$ is locally a $\P$-superdiffeomorphism affine w.r.t. the associated Levi-Civita connections.
\end{lemma}

{\it Proof:} (i) Let $\varphi_g \in \Hom\rund{sT \M, sT^* \M}$ and $\varphi_h \in \Hom\rund{sT \N, sT^* \N}$ denote the even isomorphisms defined by $\varphi_g(X) = g(\diamondsuit, X)$ and $\varphi_h(Y) = h(\diamondsuit, Y)$ for all $\P$- super vectorfields $X$ on $\M$ and $Y$ on $\N$~. Then $\varphi_g \circ (d \Phi)^\dagger \circ \rund{\Phi^* \varphi_h}$ is left-inverse to $d \Phi$~, and so $d \Phi$ is a $\P$- super vectorbundle monomorphism.

(ii) obvious by (i). $\Box$ \\

Even amoung Riemannian metrics on ordinary smooth manifolds non-trivial deformations occur, as the next example illustrates:

\begin{example}
Let $\rund{g^C}_{C \in \rz}$ be the smooth family of metrics given on $\rz^2$ by

\[
g^C_{xx} = g^C_{yy} = \exp\rund{C\rund{x^2 + y^2}} \,~, \, g^C_{xy} = 0 \,~.
\]

The Levi-Civita connection to $g^C$ is given by the Christoffel symbols

\[
{\Gamma^C}_{xx}^x = C x \,~, \, {\Gamma^C}_{xx}^y = - C y \,~, \, {\Gamma^C}_{xy}^x = C y
\]

and $x$ and $y$ interchanged, and its scalar curvature by $S^C = - 4 C$~. Now, since $g^C$ is a smooth family we may take $C \in \P_0$ and obtain a Riemannian $\P$-metric $g^C$ on $\rz^2$~. Therefore we see that, given $C, C' \in \P$~, $C \not= C'$~, there exists {\bf no} $\P$-diffeomorphism $\varphi: \rz^2 \rightarrow \rz^2$ isometric w.r.t. $g^C$ and $g^{C'}$~. In particular, any $C \in \P$ with $C^2 = 0$ yields a non-trivial infinitesimal deformation $g^C$ of the Euclidian metric on $\rz^2$~.
\end{example}

Finally, for defining the geodesic flow properly we have to be able to write $\P$- super vectorbundles on $\P$-supermanifolds as $\P$-supermanifolds themselves. Here is the general concept: To every $\varphi \in \Hom\rund{\rund{\C^\infty_\M}^{m|n}, \rund{\C^\infty_\M}^{p|q}}_0$ given by $\varphi\rund{e_k} = a_k^l f_l$~, $a_k^l \in \C^\infty(\M)_{\abs{k} + \abs{l}}$~, we associate the $\P$-supermorphism

\[
\widehat \varphi:= \rund{\Pr\nolimits_\M, \xi^k a_k^l}: \M \times \rz^{m|n} \rightarrow \M \times \rz^{p|q} \,~,
\]

$\xi^k$ denoting the standard super coordinates on $\rz^{m|n}$~. Since this assignment is again co\-va\-riant in $\varphi$ and compatible with restrictions to open subsets of $M$ we obtain a co\-va\-riant functor from $\P{\rm\bf -SuperVB}(\M)$ to the category ${\rm\bf Fam}\P{\rm\bf -SuperMan}(\M)$ of families $\Pi_\Z: \Z \twoheadrightarrow \M$ of $\P$-supermanifolds over $\M$ together with strong $\P$- super family morphisms

\[
\begin{array}{ccc}
\Z & \mathop{\longrightarrow}\limits^\Psi & \W \\
\phantom{12}_{\Pi_\Z} \twoheadsearrow & \circlearrowleft & \twoheadswarrow {}_{\Pi_\W} \phantom{1,} \\
& \M &
\end{array} \,~.
\]

Furthermore we have a 1-1-correspondence between even sections $S$ of $\E$ and sections

\[
\begin{array}{ccc}
\M & \mathop{\longrightarrow}\limits^{\widehat S} & \widehat \E \\
\phantom{12}_{\Id_\M} \searrow & \circlearrowleft & \twoheadswarrow {}_{\Pi_{\widehat \E}} \phantom{12} \\
& \M &
\end{array}
\]

of the familiy $\widehat \E$~. If $S = S^k e_k$ in a local trivialization of $\E$ then $\widehat S := \rund{\Pr\nolimits_\M, S^k}$ in the associated trivialization of $\widehat \E$~, and $\widehat{\varphi(S)} = \widehat \varphi \circ \widehat S$ for every section $S$ of $\E$ and homomorphism $\varphi: \E \rightarrow \F$ of $\P$-super vectorbundles. \\

\begin{defin} The above presented functor is called the $\P$-supermanifold realisation of $\P$- super vectorbundles on $\M$~.
\end{defin}

It obviously commutes with ${}^\#$ and ${}^\rho$~. Observe that in this realization the pullback of $\P$- super vectorbundles under a $\P$-supermorphism $\Phi: \M \rightarrow \N$ becomes precisely the pullback of families and the sum $\E \oplus \F$ of two $\P$- super vectorbundles the restricted direct product $\widehat{\E \oplus \F} = \widehat \E \times_\M \widehat \F$~. \\ % The restricted direct product $\Pi_{\Z \times_\M \W}: \Z \times_\M \W \twoheadrightarrow \M$ of two families $\Pi_\Z: \Z \twoheadrightarrow \M$ and $\Pi_\W: \W \twoheadrightarrow \M$ is defined as the sub $\P$-supermanifold of $\Z \times \W$ given by the equation $\Pi_\Z \circ \Pr_\Z = \Pi_\W \circ \Pr_\W$~, where $\Pr_\Z: \Z \times \W \twoheadrightarrow \Z$ and $\Pr_\W: \Z \times \W \twoheadrightarrow \W$ denote the canonical projections onto the first resp. the second factor, and $\Pi_{\Z \times_\M \W} := \Pi_\Z \circ \Pr_\Z = \Pi_\W \circ \Pr_\W$~. \\

In the special case $\M = \{0\}$ the $\P$-supermanifold realization gives a covariant functor from the category of the graded $\P$-modules $\P^{p|q}$~, $p, q \in \nz$~, to the category of $\P$-supermanifolds assigning to $\P^{p|q}$ the supermanifold $\rz^{p|q}$~.

\begin{theorem}[Prolongation functor]
\item[(i)] There exists a covariant functor from ${\P{\rm\bf -SuperMan}}$ to ${\rm\bf Fam}\P{\rm\bf -SuperMan}(\M)$~, called the prolongation, assigning to every $\P$-supermanifold $\M$ its super tangent bundle $\widehat{sT \M}$ and to every $\P$-supermorphism $\Phi: \M \rightarrow \N$ the $\P$-supermorphism $\widetilde \Phi: \widehat{sT \M} \rightarrow \widehat{sT \N}$ in local super charts and associated frames $\partial_i$ of $sT \M$ and $\partial_k$ of $sT \N$ defined by

\[
\widetilde \Phi := \rund{\Phi \circ \Pr\nolimits_{U^{|n}}, \xi^i \rund{\rund{\partial_i \Phi^j} \circ \Pr\nolimits_{U^{|n}} }} : U^{|n} \times \rz^{m|n} \rightarrow V^{|q} \times \rz^{p|q} \,~.
\]

\item[(ii)] Two $\P$- super vectorfields $X$ on $\M$ and $Y$ on $\N$ are $\Phi$-related iff $\widetilde \Phi \circ \widehat X = \widehat Y \circ \Phi$~.
\end{theorem}

{\it Proof:} simple calculation. $\Box$ \\

Finally also the prolongation commutes with ${}^\#$ and ${}^\rho$~. Given a $\P$-supermorphism \\
$\Phi: (\rz \times \M)|_U \rightarrow \N$~, $U \subset \rz \times M$ open, we have the choice: on the one hand we can view $\dot \Phi$ as the section $(d \Phi) \partial_t \in \rund{\Phi^* sT \M}_0(U)$~, on the other hand as the $\P$-supermorphism $\widetilde \Phi\rund{\diamondsuit, \partial_t}: (\rz \times \M)|_U \rightarrow \widehat{sT \M}$~.

\section{Integrating $\mathcal{P}$- super vectorfields}

Let $X \in \fX(\M)_0$~. If $X \in \fX(\M)$ is not purely even we do not have to give up: by the canonical embedding $\iota: \P \hookrightarrow \P \boxtimes \bigwedge \rz$ we can associate to it $X_0 + \alpha X_1 \in \fX\rund{\M^\iota}_0$~, $\alpha$ being the odd generator of $\bigwedge \rz$~, and so requiring $X$ to be even is not really a restriction thanks to the parametrization! In the following we will study in detail the differential equation 

\begin{equation} \label{diffeq}
\dot J = J^* X \,~,
\end{equation}

where $J: (\rz \times \N)_\Omega \rightarrow \M$ denotes a $\P$-supermorphism, $\Omega \subset \rz \times N$ open. Obviously if we have a solution $J$ to (\ref{diffeq}) then of course $\rund{J^\#}^\cdot = \rund{J^\#}^* X$~, and $\rund{J^\rho}^\cdot = \rund{J^\rho}^* X$~.

\begin{defin}
A $\P$-curve fulfilling (\ref{diffeq}) is called an integral curve to $X$~.
\end{defin}

First we study (\ref{diffeq}) locally, so for a moment we may assume that $\M = U^{|n}$~, $U \subset \rz^m$ open, and $\N = V^{|q}$~, $V \subset \rz^p$ open. Then we can write $X = A^i \partial_{i|} + \Delta^j \partial_{|j}$ with

\begin{eqnarray*}
&& A = \rund{A^1, \dots, A^m} \in \rund{\P \boxtimes \C^\infty\rund{U^{|n}} }_0^{\oplus m} \,~, \\
&& \Delta = \rund{\Delta^1, \dots, \Delta^n} \in \rund{\P \boxtimes \C^\infty\rund{U^{|n}} }_1^{\oplus n}
\end{eqnarray*}

and $J = \rund{f, \lambda}$ with 

\begin{eqnarray*}
&& f \in \rund{\P \boxtimes \C^\infty\rund{\Omega^{|q}} }_0^{\oplus m} = \C^\infty(\Omega)^{\oplus m} \otimes \rund{\P \boxtimes \bigwedge \rz^q}_0 \,~, \\
&& \lambda \in \rund{\P \boxtimes \C^\infty\rund{\Omega^{|q}}}_1^{\oplus n} = \C^\infty(\Omega)^{\oplus n} \otimes \rund{\P \boxtimes \bigwedge \rz^q}_1 \,~,
\end{eqnarray*}

$\Omega \subset \rz \times V$ open. Hereby $\partial_{i|}$ and $\partial_{|j}$ denote the derivatives w.r.t. the even $i$th resp. odd $j$th coordinate. Furthermore with the help of the odd coordinate functions $\xi^1, \dots \xi^n$ of $U^{|n}$ we can decompose

\begin{eqnarray*}
A = A_I \xi^I &,& A_I \in \C^\infty(U)^{\oplus m} \otimes \P_{\abs{I}} \,~, \\
\Delta = \Delta_I \xi^I &,& \Delta_I \in \C^\infty(U)^{\oplus n} \otimes \P_{1 + \abs{I}} \,~.
\end{eqnarray*}

$\P \boxtimes \bigwedge \rz^q$ is a small super algebra, so let $\n$ denote its largest ideal. Obviously $\rz \hookrightarrow \P \boxtimes \bigwedge \rz^q$ as a unital purely even subalgera, and we can decompose

\[
\P \boxtimes \bigwedge \rz^q = \rz \oplus \bigoplus_{r = 1}^\infty \I^r \,~,
\]

all $\I^r$ graded, such that $\n^r = \bigoplus_{k \geq r} \I^k$~, and so $\I^r \I^s \subset \bigoplus_{k \geq r + s} \I^k$ for all $r$ and $s \in \nz \setminus \{0\}$~. So we can continue decomposing

\begin{eqnarray*}
f = J^\# + \sum_{r = 1}^\infty f_r &,& f_r \in \C^\infty(\Omega)^{\oplus m} \otimes \I^r_0 \,~, \\
\lambda = \sum_{r = 1}^\infty \lambda_r &,& \lambda_r \in \C^\infty(\Omega)^{\oplus m} \otimes \I^r_1 \,~,
\end{eqnarray*}

and using the canonical unital even embedding $\P \hookrightarrow \P \boxtimes \bigwedge \rz^q$ also

\begin{eqnarray*}
A_I = A_I^{\#'} + \sum_{r = 1}^\infty A_{I, r} &,& A_I^{\#'} \in \C^\infty(U)^{\oplus m} \,~, \, A_{I, r} \in \C^\infty(U)^{\oplus m} \otimes \I^r_{\abs{I}} \,~, \\
\Delta_I = \Delta_I^{\#'} + \sum_{r = 1}^\infty \Delta_{I, r} &,& \Delta_I^{\#'} \in \C^\infty(U)^{\oplus m} \,~, \, \Delta_{I, r} \in \C^\infty(U)^{\oplus m} \otimes \I^r_{\abs{I} + 1} \,~.
\end{eqnarray*}

We know that $A_I^{\#'} = 0$ if $2 \not| \abs{I}$ and $\Delta_I^{\#'} = 0$ if $2 | \abs{I}$ since $\abs{A_I} = \abs{I}$ and $\abs{\Delta_I} = \abs{I} + 1$~.

Now (\ref{diffeq}) is equivalent to the ordinary order $1$ system of differential equations

\begin{eqnarray}
\dot J^\# + \sum_{r = 1}^\infty \dot f_r &=& \rund{\rund{\partial_\bk A_I^{\#'}} \circ J^\# + \sum_{s = 1}^\infty \rund{\partial_\bk A_{I, s}} \circ J^\#} \rund{\sum_{s = 1}^\infty f_s}^\bk \rund{\sum_{s = 1}^\infty \lambda_s}^I \,~, \\
\sum_{r = 1}^\infty \dot \lambda_r &=& \rund{\rund{\partial_\bk \Delta_I^{\#'}} \circ J^\# + \sum_{s = 1}^\infty \rund{\partial_\bk \Delta_{I, s}} \circ J^\#} \rund{\sum_{s = 1}^\infty f_s}^\bk \rund{\sum_{s = 1}^\infty \lambda_s}^I \,~. \notag
\end{eqnarray}

Decomposing left and right sides along $\P \boxtimes \bigwedge \rz^q = \rz \oplus \bigoplus_{r = 1}^\infty \I^r$ yields

\[
\dot J^\# = A^\# \circ J^\# \,~,
\]

which means nothing but $\dot J^\# = \rund{J^\#}^* X^\#$~, and the right hand side depends smoothly on $J^\#$~, and for all $r \in \nz \setminus \{0\}$

\begin{eqnarray*}
\dot f_r &=& \rund{\rund{\partial_i A_\emptyset^{\#'}} \circ J^\#} f_r^i + S_r\rund{J^\#, f_s, \lambda_s} \,~, \\
\dot \lambda_r &=& \rund{\Delta_{\{i\}}^{\#'} \circ J^\#} \lambda_r^i + T_r\rund{J^\#, f_s, \lambda_s} \,~,
\end{eqnarray*}

where the $S_r$~, $T_r$ are independent of $f^i_s$~, $\lambda^j_s$~, $s \geq r$~, polynomial w.r.t. $f^i_s$~, $\lambda^j_s$~, $s \leq r - 1$ and smooth w.r.t. $\rund{x, J^\#} \in V \times U$ with values in $\rund{\I^r_0}^{\oplus m}$ resp. $\rund{\I^r_1}^{\oplus n}$~. So by classical partial differential equation theory:

\begin{lemma}[Local solution of (\ref{diffeq})] \label{diffeq local}
\item[(i)] Given an open interval $I \subset \rz$~, $W \subset V$ open and $t_0 \in I$~, a solution $J: I \times W^{|q} \rightarrow U^{|n}$ of (\ref{diffeq}) is uniquely determined by $J|_{\schweif{t_0} \times W^{|q}}$~.
\item[(ii)] For each point $x_0 \in U$ there exists $\eps > 0$~, an open neighbourhood $W \subset U$ of $x_0$ and a solution $J: \ ]- \eps, \eps[ \ \times \ W^{|n} \rightarrow U^{|n}$ of (\ref{diffeq}) with $J|_{\{0\} \times W^n} = \Id_{W^{|n}}$~.
\end{lemma}

{\it Proof:} (i) by theorem 2.5.3 of \cite{Hurewicz}.

(ii) Let $x^i$~, $\xi^j$ denote the canonical super coordinate functions on $U^{|n}$~. Then we can decompose $\xi = \sum_{r = 1}^\infty \beta_r$~, $\beta_r \in \rund{\I^r_1}^{\oplus n}$~, and the initial condition $J|_{\{0\} \times W^{|n}} = \Id_{W^{|n}}$ becomes

\[
f^\#(0, \diamondsuit) = x \,~, \, f_r(0, \diamondsuit) = 0 \,~, \, \lambda_r(0, \diamondsuit) = \beta_r \,~.
\]

By theorems 2.5.7 and 2.5.9 of \cite{Hurewicz} there exists $\eps > 0$~, $W \subset U$ open and a solution \\
$f^\#: \ ]- \eps, \eps[ \ \times W \ \rightarrow U$ of the initial value problem $\dot f^\# = A_\emptyset^{\#'} \circ f^\#$ and $J^\#(0, \diamondsuit) = x$~. Since the right hand sides are linear w.r.t. $f_r(\diamondsuit, x)$ resp. $\lambda_r(\diamondsuit, x)$ and independent of all $f_s(\diamondsuit, x), \lambda_s(\diamondsuit, x)$~, $s > r$~, by the Picard-Lindelöf theorem or the example at the end of section 2.5. of \cite{Hurewicz}~, for all $x \in W$ there exist solutions $f_r(\diamondsuit, x) \in \C^\infty(] - \eps, \eps [)^{\oplus m} \otimes \I^r_0$ and $\lambda_r(\diamondsuit, x) \in \C^\infty(] - \eps, \eps [)^{\oplus n} \otimes \I^r_1$ of the initial value problem

\begin{eqnarray*}
{f_r(\diamondsuit, x)}^\cdot &=& \rund{\rund{\partial_i A_\emptyset^{\#'}} \circ f^\#}(\diamondsuit, x) f_r^i(\diamondsuit, x) + S_r\rund{f^\#(\diamondsuit, x), f_s(\diamondsuit, x), \lambda_s(\diamondsuit, x)} \,~, \\
{\lambda_r(\diamondsuit, x)}^\cdot &=& \rund{\Delta_{\{i\}}^{\#'} \circ f^\#}(\diamondsuit, x) \lambda_r^i(\diamondsuit, x) + T_r\rund{f^\#(\diamondsuit, x), f_s(\diamondsuit, x), \lambda_s(\diamondsuit, x)} \,~,
\end{eqnarray*}

$f_r(0, x) = 0$ and $\lambda_r(0, x) = \beta_r$~. By theorem 2.5.10 of \cite{Hurewicz} they depend smoothly on $x$~. $\Box$ \\

Now for the global theory we proceed as in the classical case:

\begin{lemma} \label{diffeq curve} Given an integral curve $\gamma: I \rightarrow \M$ to $X$ and $t_0 \in I^\P$~, $\gamma$ is uniquely determined by $\gamma\rund{t_0}$~.
\end{lemma}

{\it Proof:} After translation we may assume that $t_0 = 0$~. So let $\eta: I \rightarrow \M$ be a second integral curve to $X$ with $\eta(0) = \gamma(0)$~. Define

\[
U := \{t \in I \, | \, \gamma(t) = \eta(t)\} \,~.
\]

Then obviously $U$ is closed, and $0 \in U$~. But by lemma \ref{diffeq local} (i) $U$ is also open, and so $U = I$~. $\Box$ \\

Let $\Xi_X$ be the set of all $\Omega \subset \rz \times M$ open such that

\begin{itemize}
\item[\{i\}] $\Omega \cap (\rz \times \{x\}) \subset \rz$ is an interval containing $0$ for all $x \in N$ (so $\Omega$ is an interval bundle over $N$ ),
\item[\{ii\}] there exists a solution $J: (\rz \times \M)|_\Omega \rightarrow \M$ of (\ref{diffeq}) with $J|_{\{0\} \times \M} = \Id_\M$~.
\end{itemize}

\begin{theorem}[Global solution of (\ref{diffeq})]
\item[(i)] For every $\Omega \in \Xi_X$ there exists exactly one solution $J_\Omega: (\rz \times \M)|_\Omega \rightarrow \M$ of (\ref{diffeq}) with $\left.J_\Omega\right|_{\{0\} \times \M} = \Id_\M$~.
\item[(ii)] For every $x \in M$ denote by $I_{x, X^\#}$ the largest interval on which there exists an integral curve $\gamma: I \rightarrow M$ to $X^\#$ with $\gamma(0) = x$~. Then

\[
\Omega_X := \bigcup_{x \in M} I_{x, X^\#}
\]

is the largest set in $\Xi_X$~.
\end{theorem}

{\it Proof:} (i) Assume there are two solutions $J$ and $K$~. Let $x \in \M^\Q$~. Then $J(\diamondsuit, x), K(\diamondsuit, x): I \rightarrow \M$ with $I := \Omega \cap \rund{\rz \times \schweif{x^\#}}$ are integral curves to $X$~, which coincide at $0$~. So by lemma \ref{diffeq curve} they are equal. Since $x$ has been arbitrary and for large enough $\Q$ the $\Q$-points separate the $\P$-supermorphisms from $\rz \times \M$~, we have $J = K$~. \\

(ii) {\it Step I}: $\Xi_X \not= \emptyset$~. Indeed: Given $x \in M$~, by lemma \ref{diffeq local} (ii) there exist $\eps_x > 0$~, an open super coordinate neighbourhood $U_x \subset M$ of $x$ and a solution $J_x: \ ]- \eps_x, \eps_x[ \ \times \ \M|_{U_x} \rightarrow \M$ of (\ref{diffeq}) with $J|_{\{0\} \times \M|_{U_x}} = \Id_{\M|_{U_x}}$~. By lemma \ref{diffeq local} (i) all $J_x$~, $x \in M$~, coincide on their overlaps and so glue together to a solution $J: (\rz \times \M)|_\Omega \rightarrow \M$ with

\[
\Omega := \bigcup_{x \in M} \ ]- \eps_x, \eps_x[ \ \times \ U_x \,~,
\]

which so is an element of $\Xi_X$~.

{\it Step II}: $\Xi_X$ contains a largest set $\Omega_{\max}$~. Indeed: By (i) all $J_\Omega$ glue together to a solution $J_X: (\rz \times \M)|_{\Omega_{\max}} \rightarrow \M$ of (\ref{diffeq}) with $\Omega_{\max} := \bigcup_{\Omega \in \Xi_X} \Omega$~, which so is obviously the largest set in $\Xi_X$~.

{\it Step III}: $\Omega_{\max} = \Omega_X$ : `$\subset$' clear since given an arbitrary $x \in M$~, \\
$\gamma := J_X^\#\rund{\diamondsuit, x}: \rund{\rz \times \{x\}} \cap \Omega_X \rightarrow M$ is an integral curve for $X^\#$ with $\gamma(0) = x$~.

`$\supset$': For $x \in M$ write $\rund{\rz \times \{x\}} \cap \Omega_{\max} = \ ]a, b[$~, $a < 0 < b$ and assume $I_{x, X^\#} \not\subset \ ]a, b[$~. Then $a \in I_{x, X^\#}$ or $b \in I_{x, X^\#}$~. In the second case choose $\eps \in \ ]0, b[$ and an open neighbourhood $V \subset M$ of $\gamma_{x, X^\#}(b)$ such that $]- 2 \eps, 2 \eps[ \ \times V \subset \Omega_{\max}$ and $\gamma_{x, X^\#}(b - \eps) \in V$~. Finally choose an open neighbourhood $U \subset M$ of $x$ such that $J_X^\#(b - \eps, U) \subset V$~. Now define

\[
\Phi := J_X\rund{t - b + \eps, J_X\rund{b - \eps, \Pr\nolimits_\M}}: \ ]b - 3 \eps, b + \eps[ \ \times \ \M|_U \rightarrow \M \,~.
\]

Then obviously $\Phi$ fulfills the differential equation (\ref{diffeq}), and $\Phi$ and $J_X$ coincide on \\
$\{b - \eps\} \times \M|_U$~. Therefore, after maybe shrinking $U$~, by lemma \ref{diffeq local} (i) $\Phi = J_X$ on the overlap $(\rz \times \M)|_{(]b - 3 \eps, b + \eps[ \ \times \ U) \cap \Omega_{\max}}$~, which means that we can glue together $\Phi$ and $J_X$ to a solution

\[
J: (\rz \times \M)|_{\widetilde \Omega} \rightarrow \M \,~,
\]

$\widetilde \Omega := \Omega_{\max} \cup (]b - 3 \eps, b + \eps[ \ \times \ U)$~, of (\ref{diffeq}) with $J|_{\{0\} \times \M} = \Id_\M$~, so $\widetilde \Omega \in \Xi_X$ and therefore $\widetilde \Omega \subset \Omega_{\max}$~. But on the other hand $(b, x) \in \widetilde \Omega \setminus \Omega_{\max}$~. Contradiction! Similar argument in the first case $a \in I_{x, X^\#}$~. $\Box$

\begin{defin} $J_X := J_{\Omega_X}: \rund{\rz \times \M}|_{\Omega_X} \rightarrow \M$ is called the integral flow to $X$ on $\M$~, and $X$ its generator.
\end{defin}

\begin{cor}
\item[(i)] $\Omega_X = \Omega_{X^\#}$~, which so depends only on the underlying classical structure, $J_X^\# = J_{X^\#}$~, and $J_X^\rho = J_{X^\rho}$~.
\item[(ii)] Given $x \in \M^\P$~, $J_X(\diamondsuit, x): I_{x^\#, X^\#} \rightarrow \M$ is the largest integral curve to $X$ through $x$~.
\item[(iii)] $s \ \Omega_{s X^\#} = \Omega_{X^\#}$ for all $s \in \rz \setminus \{0\}$~, and $J_X\rund{s \ t, \Pr_\M} = J_{s X}$ for all $s \in \rz^\P$~.
\end{cor}

$J_X$ is really a partial action of $\rz$ on $\M$ :

\begin{theorem} \label{reallyaction} $(s + t, x) \in \Omega_{X^\#}$ for all $(s, (t, x)) \in J_{X^\#}^* \Omega_{X^\#}$ (interval bundle pullback), and on $\left.\rund{\rz^2 \times \M}\right|_{J_{X^\#}^* \Omega_{X^\#}}$

\[
J_X \circ \rund{s, J_X \circ \rund{t, \Pr\nolimits_\M}} = J_X\rund{s + t, \Pr\nolimits_\M} \,~,
\]

where $s, t: \left.\rund{\rz^2 \times \M}\right|_{\Omega_{J_X, X}} \rightarrow \rz$ denote the projections onto the first resp. second copy of $\rz$~.
\end{theorem}

This immediately implies that given $t_0 \in \rz^\P$~,

\[
\left.J_X\rund{t_0, \diamondsuit}\right|_{\M|_{\rund{\schweif{t_0} \times M} \cap \Omega_{X^\#}} }: \M|_{\rund{\schweif{t_0} \times M} \cap \Omega_{X^\#}} \rightarrow \M|_{\rund{\schweif{- t_0} \times M} \cap \Omega_{X^\#}}
\]

is a $\P$-superdiffeomorphism with inverse $\left.J_X\rund{- t_0, \diamondsuit}\right|_{\M|_{\rund{\schweif{- t_0} \times M} \cap \Omega_{X^\#}} }$~.

{\it Proof:} Let $x \in \M^\Q$ and $s, t \in \rz$ such that $\rund{s, t, x} \in J_{X^\#}^* \Omega_{X^\#}$~. Then $J_X(\diamondsuit, x)$ and \\
$J_X(\diamondsuit - t, J_X(t, x))$ are integral curves to $X$ coinciding at $t$ and so by lemma \ref{diffeq curve} on their overlaps. Therefore they glue together to an integral curve $I_{x^\#, X^\#} \cup \rund{I_{J_{X^\#}\rund{u, x^\#}} + t} \rightarrow~\M$~. We see that $s + t \in I_{J_{X^\#}\rund{s, x^\#}} + t \subset I_{x^\#, X^\#}$~, which gives the first statement, and $J_X(s, J_X(t, x)) = J_X(s + t, x)$~, which gives the second statement. $\Box$

\begin{theorem} \label{relatedintflow}
Let $Y \in \fX(\N)_0$ and $\Phi: \M \rightarrow \N$ a $\P$-supermorphism. If $X$ and $Y$ are related under $\Phi$ then $\Omega_{X^\#} \subset \rund{\Phi^\#}^* \Omega_{Y^\#}$ (interval bundle pullback), and $J_Y \circ \rund{t, \Phi \circ \Pr\nolimits_\M} = \Phi \circ J_X$ on $(\rz \times \M)|_{\Omega_{X^\#}}$~. Conversely, if $\rund{\Phi \circ J_X}^\cdot = \rund{\Phi \circ J_X}^* Y$ then $X$ and $Y$ are related under $\Phi$~.
\end{theorem}

{\it Proof:} Let $x \in \M^\Q$~. Then $\gamma := \Phi \circ J_X(\diamondsuit, x): \rund{\rz \times \schweif{x^\#}} \cap \Omega_{X^\#}$ is an integral curve to $Y$ with $\gamma(0) = \Phi(x)$ since

\[
\dot \gamma = \rund{J_X(\diamondsuit, x)^* (d \Phi)} J_X(\diamondsuit, x)^\cdot = J_X(\diamondsuit, x)^* \rund{(d \Phi) X} = J_X(\diamondsuit, x)^* \rund{\Phi^* Y} = \gamma^* Y \,~.
\]

Therefore $\rund{\rz \times \schweif{x^\#}} \cap \Omega_{X^\#} \subset I_{\Phi^\#\rund{x^\#}, Y^\#}$~, which proves the first statement. Furthermore $\gamma = J_Y\rund{\diamondsuit, \Phi(x)}$~, which proves the second statement.

For the converse just observe that

\[
(d \Phi) X = \left.(d \Phi) \dot{J_X}\right|_{\{0\} \times \M} = \left.\rund{\Phi \circ J_X}^\cdot \right|_{\{0\} \times \M} = \left.\rund{\Phi \circ J_X}^* Y \right|_{\{0\} \times \M} = \Phi^* Y \,~. \, \Box
\]

Now for $a, b \in \P_0$ such that $a^2 = b^2 = 0$ and $Y \in \fX(\M)_0$~, $J_X(a, \diamondsuit)$ and $J_Y(b, \diamondsuit)$ are $\P$-superdiffeomorphisms from $\M$ to itself with relative body $\Id_\M$~. One could wonder what is their commutator. Here the answer:

\begin{lemma} \label{commutintflow} Let $X, Y \in \fX(\M)$ and $a, b \in \P_0$ with $a^2 = b^2 = 0$~. Then
\item[(i)] $f \circ J_X(a, \diamondsuit) = f + a X f$ for all smooth super functions $f$ on $\M$~,
\item[(ii)] for all $\P$- super vectorfields $Z$ on $\M$~, $Z$ and $Z - [X,Z]$ are $J_X(a, \diamondsuit)$-related,

\item[(iii)] $\eckig{J_X(a, \diamondsuit), J_Y(b, \diamondsuit)} = J_{[Y, X]}(a b, \diamondsuit)$~.
\end{lemma}

{\it Proof:} Easy exercise in local coordinates. $\Box$ \\

In the following we will extend the theory of integral flows to arbitrary $\P$- super Lie group actions, which will lead to a super version of Palais' theorem III of section IV.2 of \cite{Palais}.

\begin{defin} A $\P$-supermanifold $\G$ with a multiplication $\P$-supermorphism \\
$m: \G \times \G \rightarrow~\G$~, an inversion $\P$-supermorphism ${}^{- 1}: \G \rightarrow \G$ and a neutral element $1 \in \M^\P$ fulfilling the usual group axioms is called a $\P$- super Lie group. The sub $\P$- super Lie algebra of right-invariant $\P$- super vectorfields $X \in \fX(\G)$ is called the super Lie algebra of $\G$~.
\end{defin}

Here $X \in \fX(\G)$ right-invariant means $\rund{d_1 m} X = m^* X$~, where $d_1$ denotes the differential w.r.t. the first entry, and is, for $\Q$ large enough, equivalent to $X$ being related to itself under all right-translations $m(\diamondsuit, g)$~, $g \in \G^\Q$~.

\begin{example} \label{superLie}
Let $\alpha \in \P_1$~. Then $\rz^{0|1}$ with multiplication $m := \xi + \eta + \xi \eta \alpha$~, inversion $- \xi$ and neutral element $0$ is a $\P$- super Lie group. Its super Lie algebra is generated by $X := (1 + \xi \alpha) \partial_\xi$~, and $[X, X] = 2 \alpha X$~. Its relative body is $\rund{\rz^{0|1}, +}$~. For $\alpha \not= 0$ there exists no $\P$- super Lie group isomorphism from $\rund{\rz^{0|1}, +}$ to $\rund{\rz^{0|1}, m}$ since the first one is abelian, the second one is not. So we have constructed a non-trivial local deformation of the ordinary super Lie group $\rund{\rz^{0|1}, +}$~.
\end{example}

From now on let $\G$ be a $\P$- super Lie group with super Lie algebra $\g$~. As in the classical case one has:

\begin{lemma} \label{Liealgebra}
\item[(i)] There is a canonical $\P$-linear isomorphism $\g \simeq sT_1 \G$ given by $X \mapsto X(1)$ and $\rund{d_1 m}(1, \diamondsuit) v \mapsfrom v$~.
\item[(ii)] There is a unique map $\exp: \widehat \g \rightarrow \G$ such that $\exp(0) = 1$~, $(d \exp)(0) = \Id_\g$ and $\exp\rund{(s + t) \Pr\nolimits_{\widehat \g}} = \exp\rund{s \Pr\nolimits_{\widehat \g}} \exp\rund{t \Pr\nolimits_{\widehat \g}}$ on $\rz^2 \times \widehat \g$~.
\item[(iii)] Let $X, Y \in \g_0$ and $a, b \in \P_0$ such that $a^2 = b^2 = 0$~. Then $J_X = m \circ \rund{\exp(t X), \Pr\nolimits_\G}$~, and $\eckig{\exp(a X), \exp(b Y)} = \exp(a b [Y, X])$~.
\end{lemma}

{\it Proof:} (i) obvious.

(iii) By the properties of $\exp$ and since $X$ is right-invariant

\[
\dot \Phi = \left.\partial_u\rund{\exp(u X) \Phi}\right|_{u = 0} = \Phi^* \rund{\rund{d_1 m}(1, \diamondsuit) X} = \Phi^* X \,~,
\]

which proves the first statement. The second follows directly from the first and lemma \ref{commutintflow}.

(ii) {\it Uniqueness}: obvious by the first statement in (iii) by passing from $\P$ to $\Q$ large enough and using the isomorphism $\g \simeq \rund{\rund{\bigwedge \rz} \boxtimes \g}_0$ of graded vectorspaces.

{\it Existence}: Take a graded base $\rund{\xi_i}$ of $\g$ and associated super coordinates $c^i$ on $\widehat \g$~. Let \\
$R := c^i \otimes \xi_i \in \fX\rund{\widehat \g \times \G}$~. $X = X^i \xi_i \in \rund{\Q \boxtimes_\P \g}_0$ given, $R$ and $X$ are related under the embedding $\rund{\widehat X, \Id_\G}: \G \hookrightarrow \widehat \g \times \G$~, and $X^\#$ is complete since $J_{X^\#} = \rund{\exp\rund{t X^\#}, \Pr\nolimits_G}$ is its integral flow. Therefore $\rz \times G = \Omega_{X^\#} \subset \rund{\rz \times \schweif{X^\#} \times G} \cap \Omega_{R^\#}$~, and \\
$J_R \circ \rund{t, X, \Pr\nolimits_\G} = \rund{X, J_{\rho(X)} }$ on $\rz \times \G$~. Since moreover $X$ has been arbitrary we see that in particular $R^\#$ is complete. Define

\[
\exp := \Pr\nolimits_\G \circ J_R \circ \rund{1, \Pr\nolimits_{\widehat \g}, 1} \,~.
\]

Since $m\rund{\exp(t X), \Pr\nolimits_\G} = J_X$ for all $X \in \g$~, the desired properties of $\exp$ follow immediately. $\Box$ \\

Given a connected $\P$- super Lie group $\G$ of super dimension $(1, 0)$ with $\P$- super Lie algebra $\g$~, $\exp(t, X)$ gives an isomorphism from either $(\rz, +)$ or $\rund{\rz / \zz, +}$ to $\G$ after appropriate choice of the generator $X$ of $\g$~. So, in contrast to example \ref{superLie}, super Lie groups of dimension $(1, 0)$ admit no non-trivial local deformations.

\begin{prop}
Let $\Phi: \G \times \M \rightarrow \M$ be a $\P$- super action. Then

\[
\varphi:= \rund{d_1 \Phi}\rund{1, \Pr\nolimits_\M}: \g \rightarrow \fX(\M)
\]

is an even $\P$- super Lie algebra homomorphism, and $\Phi\rund{\exp(t X), \Pr\nolimits_\M} = J_{\varphi(X)}$ for all $X \in \g$~.
\end{prop}

{\it Proof:} Since $\Phi$ is a $\P$- super action, $\Phi(1, \diamondsuit) = \Id_\M$ and

\begin{eqnarray*}
\rund{\Phi\rund{\exp(t X), \Pr\nolimits_\M}}^\cdot &=& \left.\partial_u \rund{\Phi \circ \rund{\exp(u X), \Phi\rund{\exp(t X), \Pr\nolimits_\M}} }\right|_{u = 0} \\
&=& \Phi^* \rund{\left.\partial_u \Phi\rund{\exp(u X), \Pr\nolimits_\M}\right|_{u = 0}} \\
&=& \Phi^* \varphi(X) \,~,
\end{eqnarray*}

which proves the second statement. Therefore after passing from $\P$ to $\Q$ large enough for all $X, Y \in \rund{\Q \boxtimes_\P \g}_0$ and $a, b \in \Q_0$ with $a^2 = b^2 = 0$ by lemma \ref{commutintflow}

\begin{eqnarray*}
J_{[\varphi(Y), \varphi(X)]}(a b, \diamondsuit) &=& \eckig{J_{\varphi(X)}(a, \diamondsuit), J_{\varphi(Y)}(b, \diamondsuit)} \\
&=& \eckig{\Phi\rund{\exp(a X), \Pr\nolimits_\M}, \Phi\rund{\exp(b Y), \Pr\nolimits_\M}} \\
&=& \Phi\rund{\eckig{\exp(a X), \exp(b Y)}, \Pr\nolimits_\M} \\
&=& \Phi\rund{\exp(a b [Y, X]), \Pr\nolimits_\M} \,~.
\end{eqnarray*}

Since in a local super chart the left hand side is equal to $\Id_\M + a b [\varphi(Y), \varphi(X)]$ and the right hand side to $\Id_\M + a b \varphi([Y, X])$ we have $[\varphi(Y), \varphi(X)] = \varphi([Y, X])$ for all $X, Y \in \rund{\Q \boxtimes_\P \g}_0$ and so for all $X, Y \in \g$~. $\Box$ \\

The converse is also true:

\begin{theorem}[Super Palais theorem] \label{Palais}
Let $\G$ be a simply connected $\P$- super Lie group with super Lie algebra $\g$ and $\varphi: \g \rightarrow \fX(\M)$ an even $\P$- super Lie algebra homomorphism such that $\varphi(\g)^\# \subset \fX(M)$ consists of complete vectorfields. Then there exists a unique $\P$- super action $J_\varphi: \G \times \M \rightarrow \M$ such that $\left.\rund{d J_\varphi}\right|_{\{1\} \times \M} X = \varphi(X)$ for all $X \in \g$~.
\end{theorem}

Moreover, theorem III of section IV.2 of \cite{Palais} says that infact are equivalent

\begin{itemize}
\item[(i)] $\varphi(\g)^\#$ consists of complete vectorfields,
\item[(ii)] there exists a basis $\rund{X_i}$ of $g^\#$ such that all $\varphi^\#\rund{X_i}$ are complete.
\end{itemize}

The proof of theorem \ref{Palais} needs some more preparation:

\begin{lemma} \label{vectorfieldstoMP}
\item[(i)] There exists a unique $\rz$-linear sheaf embedding ${}^\P: (sT \M)_0 \rightarrow \rund{\Pi_{\M^\P}}_* T \M^\P$ such that for every $U \subset M$ open, $X \in \fX\rund{\M|_U}$ and integral curve $\gamma: I \rightarrow \M_U$ to $X$~, $\left.\gamma^\P\right|_I$ is an integral curve to $X^\P$~.
\item[(ii)] Given $x \in \M^\P$~, ${}^\P$ induces an isomorphism ${}^\P: \rund{sT_x \M}_0 \mathop{\rightarrow}\limits^\sim T_x \M^\P$ such that

\[
\begin{array}{ccc}
\phantom{12345} \rund{(sT \M)_0}_{x^\#} & \mathop{\longrightarrow}\limits^{{}^\P} & \rund{T \M^\P}_{x^\#} \phantom{123} \\
\diamondsuit(x) \twoheaddownarrow & \circlearrowleft & \twoheaddownarrow \diamondsuit(x) \\
\phantom{12345} \rund{sT_x \M}_0 & \mathop{\mathop{\longrightarrow}\limits_{\sim}}\limits_{{}^\P} & T_x \M^\P \phantom{1234}
\end{array} \,~.
\]

${\dot \gamma(0)}^\P = \rund{\left.\gamma^\P\right|_I}^\cdot(0)$ for every $\P$-curve $\gamma: I \rightarrow \M$~.
\item[(iii)] If $X \in \fX(\M)$ and $Y \in \fX(\N)$ are related under the $\P$-supermorphism $\Phi: \M \rightarrow \N$ then so are $X^\P$ and $Y^\P$ under $\Phi^\P$~.
\item[(iv)] For every $X \in \fX(\M)$

\[
\left.J_X^\P\right|_{\rz \times \M^\P} = J_{X^\P} \,~.
\]

\item[(v)] ${}^\P$ is a Lie algebra sheaf embedding.
\end{lemma}

{\it Proof:} Take a local super chart $U^{m|n}$~, $U \subset \rz^m$ open, of $\M$ and bases $\rund{b_d}$ and $\rund{c_e}$ of $\m_0$ resp. $\m_1 = \P_1$~. Let $V \subset U$ be open and $X = A^i \partial_{i|} + \Delta^j \partial_{|j} \in \rund{\P \boxtimes \fX\rund{V^{|n}}}_0$~, $A^i \in \rund{\P \boxtimes \C^\infty\rund{V^{|n}}}_0$~, $\Delta^j \in \rund{\P \boxtimes \C^\infty\rund{V^{|n}}}_1$~. Then we can decompose

\[
\rund{A^i}^\P = \rund{A^i}^\# \circ {}^\# + A^{i, d} b_d \phantom{12}~, \phantom{12} \rund{\Delta^j}^\P = \Delta^{j, e} c_e \phantom{12}~, \phantom{12} A^{i, d}~, \Delta^{j, e} \in \C^\infty\rund{\rund{V^{|n}}^\P} \,~.
\]

Let $\partial_i$~, $\partial_{i|, d}$ and $\partial_{|j, e}$ denote the partial derviatives on $\rund{U^{|n}}^\P = U \times \rund{\m_0^{\oplus m} \oplus \m_1^{\oplus n}}$ w.r.t. the $i$-th coordinate on $U$~, the $b_d$-coefficient in the $i$-th entry of $\m_0^{\oplus m}$ and the $c_e$-coefficient in the $j$-th entry of $\m_1^{\oplus n}$ resp.. Then a simple calculation shows that (i) is fulfilled iff

\[
X^\P := \rund{\rund{A^i}^\# \circ {}^\#} \partial_i + A^{i, d} \partial_{i|, d} + \Delta^{j, e} \partial_{|j, e} \,~.
\]

Also (ii) is fulfilled by an easy calculation. \\

Now let $X \in \rund{\P \boxtimes \fX\rund{U^{|n}} }_0$ and $Y \in \rund{\P \boxtimes \fX\rund{V^{|q}} }_0$ be related under the $\P$-supermorphism $\Phi: U^{|n} \rightarrow V^{|q}$~, $x \in \M^\P$ and $\gamma: I \rightarrow \M$ an integral curve to $X$ with $\gamma(0) = x$~. Then by theorem \ref{relatedintflow} $\Phi \circ \gamma$ is an integral curve to $Y$~. So by (i) $\left.\gamma^\P\right|_I$ is an integral curve to $X^\P$ and $\left.\Phi^\P \circ \gamma^\P\right|_I = \left.\rund{\Phi \circ \gamma}^\P\right|_I$ is an integral curve to $Y^\P$~. Therefore

\[
Y^\P(\Phi(x)) = \rund{\partial_t \rund{\left.(\Phi \circ \gamma)^\P\right|_I}}(0) = \rund{\partial_t \rund{\left.\Phi^\P \circ \gamma\right|_I}}(0) = \rund{\rund{d \Phi^\P} X^\P}(x) \,~,
\]

and so $X^\P$ and $Y^\P$ are related under $\Phi^\P$~. This shows that ${}^\P$ is infact globally defined and fulfills (iii). Furthermore let $x \in \M^\P$~. Then $\left.J_X^\P\right|_{\rz \times \M^\P}(\diamondsuit, x)$ is an integral curve to $X^\P$ and so must equal $J_{X^\P}(\diamondsuit, x)$~. This proves (iv). Finally for proving (v) we obtain by lemma \ref{commutintflow} that for all $x \in \M^\P$ and $t, u \in \rz$ sufficiently small in a local super chart of $\M$

\begin{eqnarray*}
J_{[Y, X]}(t u, \diamondsuit)^\P(x) &=& \eckig{J_X(t, \diamondsuit), J_Y(u, \diamondsuit)}^\P(x) + O\rund{\norm{(t, u)}^3} \\
&=& \eckig{J_X(t, \diamondsuit)^\P, J_Y(u, \diamondsuit)^\P}(x) + O\rund{\norm{(t, u)}^3} \\
&=& J_{\eckig{Y^\P, X^\P}}(t u, x) + O\rund{\norm{(t, u)}^3} \,~,
\end{eqnarray*}

and so $[Y, X]^\P = \eckig{Y^\P, X^\P}$~. $\Box$ \\

{\it Proof of theorem \ref{Palais}:} By lemmas \ref{Liealgebra} (i) and \ref{vectorfieldstoMP} (ii) and (iv) the Lie algebra of $\G^\Q$ is given by $\rund{\Q \boxtimes_\P \g}_0$~, and by lemma \ref{vectorfieldstoMP} (v) we have a Lie algebra homomorphism

\[
\rund{\Q \boxtimes_\P \g}_0 \rightarrow \fX\rund{\M^\Q} \,~, \, X \mapsto \varphi(X)^\Q \,~,
\]

where we have already extended $\varphi$ to an even $\Q$- super Lie algebra homomorphism $\Q \boxtimes_\P \g \rightarrow \Q \boxtimes_\P \fX(\M)$~. Let $X \in \rund{\Q \boxtimes_\P \g}_0$ be arbitrary. Then since $\varphi(X)^\#$ is complete, we have an integral flow $J_{\varphi(X)}: \rz \times \M \rightarrow \M$~, which is a $\Q$-supermorphism. $\left.J_{\varphi(X)}^\Q\right|_{\rz \times \M^\Q}: \rz \times \M^\Q \rightarrow \M^\Q$ is the integral flow to $\varphi(X)^\Q$~, which shows that also $\varphi(X)^\Q$ is complete. Therefore and since with $\G$ also $\G^\Q$ is simply connected by lemma \ref{homot}, by Palais' theorem, theorem III of section IV.2 of \cite{Palais}, there exists a unique smooth Lie group action $\sigma: \G^\Q \times \M^\Q \rightarrow \M^\Q$ such that $(d \sigma)(1, x) X = \varphi(X)(x)$ for all $x \in \M^\Q$ and $X \in \rund{\Q \boxtimes_\P \g}_0$~. Moreover $J_{\varphi(X)}^\Q(t, x) = \sigma\rund{\exp_{\G^\Q}(t X), x}$ for all $x \in \M^\Q$~, $X \in \rund{\Q \boxtimes_\P \g}_0$ and $t \in \rz$~. \\

{\it Uniqueness:} $J_\varphi^\Q: \G^\Q \times \M^\Q \rightarrow \M^\Q$ is a smooth action such that for all $X \in \rund{\Q \times_\P \g}_0$

\[
\left.\rund{d J_\varphi^\Q} X^\Q\right|_{\{1\} \times \M^\Q} X = \varphi(X)^\Q~,
\]

and so must equal $\sigma$~. \\

{\it Existence:} Take a graded base $\rund{\xi_i}$ of $\g$ and associated super coordinates $c^i$ on $\widehat \g$~. Define $R := c^i \otimes \varphi\rund{\xi_i} \in \fX\rund{\widehat \g \times \M}$~. As in the proof of the existence part of lemma \ref{Liealgebra} (ii) one deduces from the completeness of all $\varphi(X)^\#$ the completeness of $R^\#$ and shows that $J_R \circ \rund{t, \widehat X, \Pr\nolimits_\M} = \rund{\widehat X, J_{\varphi(X)} }$ on $\rz \times \widehat{s T \M}$ for all $X \in \rund{\Q \boxtimes_\P \g}_0$~. \\

Let $V_0 \subset \g^\#$ and $U_0 \subset G$ be open neighbourhoods of $0$ resp. $1$ such that $\exp_\G: \left.{\widehat \g}\right|_{V_0} \rightarrow \G|_{U_0}$ is a $\P$-superdiffeomorphism. Define

\[
\Phi_0 := \Pr\nolimits_\M \circ J_R \circ \rund{1, \exp_\G^{- 1} \circ \Pr\nolimits_\G, \Pr\nolimits_\M}: \G|_{U_0} \times \M \rightarrow \M \,~.
\]

Then $\Phi_0^\Q$ is a restriction of $\sigma$~, indeed: For all $X \in {\widehat \g}^\Q$ such that $X^\# \in V_0$

\[
\Phi_0\rund{\exp_\G^\Q X, \diamondsuit} = J_{\varphi(X)}(1, \diamondsuit)
\]

as $\Q$-supermorphisms $\M \rightarrow \M$~. In particular for all $x \in \M^\Q$

\[
\Phi_0^\Q\rund{\exp_{\G^\Q} X, x} = \Phi^\Q\rund{\exp_{\G}^\Q X, x} = J_{\varphi(X)}^\Q(1, x) = \sigma\rund{\exp_{\G^\Q} X, x} \,~.
\]

All $\P$-supermorphisms $\Phi: \G|_U \times \M \rightarrow \M$~, $U \subset G$ open, such that $\Phi^\Q$ is a restriction of $\sigma$ coincide on their overlaps, and so they glue together to a largest $\P$-supermorphism $J_\varphi: \G|_{U_{\max}} \times \M \rightarrow \M$~, $U_{\max} \subset G$ open, such that $J_\varphi^\Q$ is also a restriction of $\sigma$~. At least $1 \in U_0 \subset U_{\max}$~, so $U_{\max} \not= \emptyset$~. But $U_{\max}$ is also closed:

\begin{quote}
Assume $a \in \G^\P$ such that $a^\# \in \partial U_{\max}$ and choose $b \in \G^\P$ such that \\
$b^\# \in U_{\max} \cap \rund{U_0^{- 1} a^\#}$ and so $a^\# \in U_0 b^\#$~. Define the $\P$-supermorphism

\[
\Psi := \Phi_0 \circ \rund{\Pr\nolimits_\G b^{- 1}, J_\varphi \circ \rund{b, \Pr\nolimits_\M}}: \G|_{U_0 b^\#} \times \M \rightarrow \M \,~.
\]

Then $\Psi^\Q$ is a restriction of $\sigma$ since for all $j \in \G^\Q$ and $x \in \M^\Q$

\[
\Psi^\Q(j, x) = \Phi_0^\Q\rund{j b^{- 1}, J_\varphi(b, x)} = \sigma\rund{j b^{- 1}, \sigma(b, x)} = \sigma(j, x) \,~.
\]

Therefore $a^\# \in U_0 b^\# \subset U_{\max}$~.
\end{quote}

So we have $U_{\max} = G$~. $J_\varphi$ is indeed an action of $\G$ since for all $g, h \in \G^\Q$ and $x \in \M^\Q$

\[
J_\varphi^\Q(g, J_\varphi(h, x)) = \sigma(g, \sigma(h, x)) = \sigma(g h, x) = J_\varphi^\Q(g h, x) \,~.
\]

Finally let $X \in \rund{\Q \boxtimes_\P \g}_0$~. Then for small $t \in \rz$

\[
\rund{J_\varphi} \circ \rund{\exp(t X), \Pr\nolimits_\M} = \Phi_0\rund{\exp_\G^\P(t X), \diamondsuit} = J_{\varphi(t X)}(1, \diamondsuit) = J_{\varphi(X)}(t, \diamondsuit) \,~,
\]

and so after taking $\Q$ large enough we obtain $\left.d \rund{J_\varphi} X\right|_{\{1\} \times \M} = \varphi(X)$ for all $X \in \g$~. $\Box$ \\

This result explains and generalizes the work of J. Monterde and O. A. Sánchez-Valenzuela of \cite{MontSan} in the complete case, who did this for three Lie group structures on $\G := \rz^{1|1}$~. \\

For the rest of this section let $\M$ and $\N$ be equipped with affine $\P$-connections $\nabla$ and $\nabla^\N$ resp. Riemannian $\P$-supermetrics $g$ and $h$~. 

\begin{defin}
Let $\B$ be a $\P$-supermanifold with body $B$~, $\Omega \subset B \times M$ open and \\
$\Phi: (\B \times \M)|_\Omega \rightarrow \N$ a $\P$-supermorphism, so a family of partial $\P$-supermorphisms from $\M$ to $\N$~. Then $\Phi$ is called
\item[(i)] affine w.r.t. $\nabla$ and $\nabla^\N$ iff for all $\P$- super vectorfields $Y, Z$ on $\M$

\[
\rund{d J_X} \rund{\nabla_Y Z} = \rund{J_X^* \nabla^\N}_Y \rund{\rund{d J_X} Z} \,~,
\]

where we used the canonical embedding $sT \M \hookrightarrow \Pr_\M^* sT \M$~, $\Pr_\M: \B \times \M \twoheadrightarrow \M$ denoting the canonical projection, as sections constant in $\B$-direction,

\item[(ii)] isometric w.r.t. $g$ and $h$ iff for all $\P$- super vectorfields $Y, Z$ on $\M$

\[
\rund{\Phi^* h} \rund{(d \Phi) Y, (d \Phi) Z} = g(Y, Z) \,~.
\]

\end{defin}

The following lemma is obvious:

\begin{lemma}
Let $\Phi: (\B \times \M)|_\Omega$ a $\P$-supermorphism, so a family of partial $\P$-supermorphisms from $\M$ to $\N$~, and let $\M$ and $\N$ be equipped with affine $\P$-connections resp. Riemannian $\P$-supermetrics. Let $\Q$ be a small algebra. Then:

If $\Phi$ is affine resp. isometric then all $\Phi(b, \diamondsuit): \M|_{\rund{\schweif{b^\#} \times M} \cap \Omega} \rightarrow \N$~, $b \in \B^\Q$~, (they are $(\Q \boxtimes \P)$-supermorphisms) are affine resp. isometric. If the $\Q$-points of $\B$ separate the $\P$-functions on $\B$ then also the converse is true.
\end{lemma}

Now one can ask when $J_X$~, $X \in \fX(\M)_0$~, is affine resp. isometric. Here is the answer:

\begin{theorem} \label{Killingisom} Let $X \in \fX(\M)_0$~. Then $J_X$ is affine resp. isometric iff for all $\P$- super vectorfields $Y, Z$ on $\M$

\[
\eckig{X, \nabla_Y Z} = \nabla_{[X, Y]} Z + \nabla_Y [X, Z] \, \text{ resp. } \, X g(Y, Z) = g([X, Y], Z) + g(Y, [X, Z]) \,~.
\]

\end{theorem}

{\it Proof:} Let $U \subset M$ be open and $Y, Z \in \fX\rund{\M|_U}$~. Take $t_0 \in \rz$ and $a \in \P_0$ such that $a^2 = 0$~. Define $\Phi := J_X\rund{t_0, \diamondsuit}$ and

\[
W := \rund{\Phi^{- 1}}^* ((d \Phi) Y) \, \text{ and } \, R := \rund{\Phi^{- 1}}^* ((d \Phi) Z) \in \fX\rund{\M|_{\rund{\{t_0\} \times U} \cap \Omega_{X^\#} }} \,~.
\]

Then by lemma \ref{commutintflow} (ii) $Y$ and $W - a [X, W]$ and $Z$ and $R - a [X, R]$ are $J_X\rund{t_0 + a, \diamondsuit}$-related. \\

\begin{quote}
For the Riemannian case define

\[
f := \rund{J_X^* g} \rund{\rund{d J_X} Y, \rund{d J_X} Z} \in \C^\infty_{\rz \times \M}\rund{(\rz \times U) \cap \Omega_{X^\#}} \,~.
\]

By lemma \ref{commutintflow} (i)

\begin{eqnarray*}
f\rund{t_0 + a, \diamondsuit} &=& \rund{J_X\rund{t_0 + a, \diamondsuit}^* g} \rund{\rund{d J_X\rund{t_0 + a, \diamondsuit}} Y, \rund{d J_X\rund{t_0 + a, \diamondsuit}} Z} \\
%&=& \rund{J_X\rund{t_0 + a, \diamondsuit}^* g} \rund{J_X\rund{t_0 + a, \diamondsuit}^* (W - a [X, W]), J_X\rund{t_0 + a, \diamondsuit}^* (R - a [X, R])} \\
&=& g(W - a [X, W], R - a [X, R]) \circ J_X\rund{t_0 + a, \diamondsuit} \\
%&=& \rund{g(W, R) \circ J_X\rund{a, \diamondsuit} - a g(W, [X, R]) - a g([X, W], R)} \circ J_X\rund{t_0, \diamondsuit}  \\
%&=& g(W, R) \circ J_X\rund{t_0, \diamondsuit} + a \rund{X g(W, R) - g(W, [X, R]) \circ J_X\rund{t_0, \diamondsuit} - g([X, W], R) \circ J_X\rund{t_0, \diamondsuit}} \\
&=& f\rund{t_0, \diamondsuit} + a \rund{X g(W, R) - g([X, W], R) - g(W, [X, R])} \circ \Phi \,~.
\end{eqnarray*}

Without restriction we may assume that there exists $a \in \P_0 \setminus \{0\}$ such that $a^2 = 0$~, and so we obtain

\[
\dot f\rund{t_0, \diamondsuit} = \rund{X g(W, R) - g([X, W], R) - g(W, [X, R])} \circ \Phi \,~.
\]

For the affine case define

\[
f := \rund{d_\M J_X}^{- 1} \rund{J_X^* \nabla}_Y \rund{\rund{d J_X} Z} \in \rund{\Pr\nolimits_\M^* sT \M} \rund{(\rz \times U) \cap \Omega_{X^\#}} \,~,
\]

where $\Pr_\M: \rz \times \M \twoheadrightarrow \M$ denotes the canonical projection and $d_\M$ the super differential restricted to $\Pr^* sT \M$~. By using lemma \ref{commutintflow} (ii) one easily checks that

\[
(d \Phi)^{- 1} \Phi^* \rund{\nabla_W R + a \rund{\eckig{X, \nabla_W R} - \nabla_{[X, W]} R - \nabla_W [X, R]} }
\]

and $\nabla_W R - a \rund{\nabla_{[X, W]} R + \nabla_W [X, R]}$ are $J_X\rund{t_0 + a, \diamondsuit}$-related. Therefore using lemma \ref{pullbackconnright} now

\begin{eqnarray*}
f\rund{t_0 + a, \diamondsuit} &=& \rund{d J_X\rund{t_0 + a, \diamondsuit}}^{- 1} \rund{J_X\rund{t_0 + a, \diamondsuit}^* \nabla}_Y \rund{\rund{d J_X\rund{t_0 + a, \diamondsuit}} Z} \\
&=& \rund{d J_X\rund{t_0 + a, \diamondsuit}}^{- 1} J_X\rund{t_0 + a, \diamondsuit}^* \rund{\nabla_{W - a [X, W]} (R - a [X, R]) } \\
&=& f\rund{t_0, \diamondsuit} + a (d \Phi)^{- 1} \Phi^* \rund{\eckig{X, \nabla_W R} - \nabla_{[X, W]} R - \nabla_W [X, R]} \,~,
\end{eqnarray*}

and so

\[
\dot f\rund{t_0, \diamondsuit} = (d \Phi)^{- 1} \Phi^* \rund{\eckig{X, \nabla_W R} - \nabla_{[X, W]} R - \nabla_W [X, R]} \,~.
\]

\end{quote}

For $t_0 = 0$ we obtain $\Phi = \Id_\M$ and so $W = Y$ and $R = Z$~, which implies $f(0, \diamondsuit) = g(Y, Z)$ in the Riemannian and $f(0, \diamondsuit) = \nabla_Y Z$ in the affine case. So since $(\rz \times U) \cap \Omega_{X^\#}$ is an interval bundle containing $\{0\} \times U$ the statement has become obvious. $\Box$ \\

This theorem motivates the following definition:

\begin{defin} Let $X$ be a homogeneous $\P$- super vectorfield on $\M$~.
\item[(i)] $X$ is called infinitesimally affine iff for all $\P$- super vectorfields $Y, Z$ on $\M$~, $Y$ homogeneous,

\[
\eckig{X, \nabla_Y Z} = \nabla_{[X, Y]} Z + (- 1)^{\abs{X} \abs{Y}} \nabla_Y [X, Z] \,~.
\]

\item[(ii)] $X$ is called Killing iff for all $\P$- super vectorfields $Y, Z$ on $\M$~, $Y$ homogeneous,

\[
X g(Y, Z) = g([X, Y], Z) + (- 1)^{\abs{X} \abs{Y}} g(Y, [X, Z])
\]

or equivalently using the Levi-Civita connection $\nabla$ to $g$

\[
g\rund{\nabla_Y X, Z} + (- 1)^{\abs{X} \abs{Z}} g\rund{Y, \nabla_Z X} = 0 \,~.
\]

\item[(iii)] An arbitrary $\P$- super vectorfield on $\M$ is called infinitesimally affine resp. Killing iff so are its even and odd part. We denote the graded $\P$-submodule of all infinitesimally affine resp. Killing $\P$- super vectorfields by $\fX_\infaff(\M)$ resp. $\fX_\Kill(\M)$~.
\end{defin}

\begin{cor}
\item[(i)] $\fX_\infaff(\M) \subset \fX(\M)$ is a sub $\P$- super Lie algebra.
\item[(ii)] If $\M$ is a $\P$- Riemannian supermanifold equipped with the Levi-Civita connection then $\fX_\Kill(\M) \subset \fX_\infaff(\M)$ is a sub $\P$- super Lie algebra.
\end{cor}

{\it Proof:} Direct consequence of theorems \ref{Killingisom}, lemma \ref{commutintflow}, and lemma \ref{isometric} (ii). $\Box$

\begin{cor}
Let $\varphi: \g \rightarrow \fX(\M)$ and $J_\varphi$ as in theorem \ref{Palais}. Then $J_\varphi$ is affine resp. isometric iff $\varphi(\g) \subset \fX_{\infaff}(\M)$ resp. $\varphi(\g) \subset \fX_{\Kill}(\M)$~.
\end{cor}

{\it Proof:} `$\Rightarrow$': obvious by theorem \ref{Killingisom}~.

`$\Leftarrow$': As in the proof of theorem \ref{Palais} let $V_0 \subset {\widehat \g}^\#$ and $U_0 \subset G$ be open neighbourhoods of $0$ resp. $1$ such that $\exp_\G: \left.{\widehat \g}\right|_{V_0} \rightarrow \G|_{U_0}$ is a $\P$-superdiffeomorphism. $J_\varphi(\exp X) = J_{\varphi(X)}(1, \diamondsuit)$~, which is affine resp. an isometry by theorem \ref{Killingisom}, for all $X \in (\Q \boxtimes \g)_0 = \rund{\widehat g}^\Q$~. Therefore $\left.J_\varphi\right|_{\G|_{U_0} \times \M}$ is affine resp. isometric.

Now let $U \subset G$ denote the largest open subset such that $\left.J_\varphi\right|_{\G|_U \times \M}$ is affine resp. isometric. At least $U_0 \subset U$~, which is therefore non-empty. As in the proof of theorem \ref{Palais} we see that $U$ is also closed:

\begin{quote}
Assume $a \in \G^\P$ such that $a^\# \in \partial U$ and take $b \in \G^\P$ such that $b^\# \in U \cap \rund{U_0^{- 1} a^\#}$ and so $a^\# \in U_0 b^\#$~. For all $c \in \rund{\G|_{U_0 b^\#}}^\Q$ we have

\[
J_\varphi(c, \diamondsuit) = J_\varphi\rund{c b^{- 1}, \diamondsuit} \circ J_\varphi(b, \diamondsuit) \,~,
\]

which is affine resp. an isometry since $c b^{- 1} \in U_0$ and $b \in U$~. Therefore $\left.J_\varphi\right|_{\G|_{U_0 b^\#} \times \M}$ is affine resp. isometric, so $a^\# \in U$~.
\end{quote}

Since $G$ is connected we have $U = G$~. $\Box$

\section{The geodesic flow on a $\mathcal{P}$-supermanifold with affine $\mathcal{P}$-connection}

Throughout this section let $\nabla$ be an affine $\P$-connection on $\M$~.

\begin{defin} \label{geodflow}
Let $U \in \fX\rund{\widehat{sT \M}}_0$ be given in local super coordinates by

\[
U := \xi^i \partial_{x^i} - \xi^j \xi^i \rund{\Gamma_{i j}^k \circ \Pi_{\widehat{sT \M}} } \partial_{\xi^k}~,
\]

where $\Gamma_{i j}^k$ denote the Christoffel symbols of $\nabla$~. Then the integral flow

\[
\Phi := J_U : \left.\rund{\rz \times \widehat{sT \M}}\right|_\Omega \rightarrow \widehat{sT \M}
\]

with $\Omega := \Omega_{U^\#}$ is called the geodesic flow on $\M$ w.r.t. $\nabla$~.
\end{defin}

$U$ and so also $\Phi$ are indeed globally defined:

\begin{lemma} \label{geodcompatibleaffinelocal}
Let $\Xi: V^{|n} \rightarrow W^{|q}$~, $V \subset \rz^m$~, $W \subset \rz^p$ open, be a $\P$-supermorphism, $\nabla^V$ and $\nabla^W$ affine $\P$-connections on $V^{|n}$ res. $W^{|q}$ given by the Christoffel symbols \\
$\Gamma^{V, k}_{i j} \in \rund{\P \boxtimes \C^\infty\rund{V^{|n}} }_{\abs{i} + \abs{j} + \abs{k}}$ resp. $\Gamma^{W, k}_{i j} \in \rund{\P \boxtimes \C^\infty\rund{W^{|q}} }_{\abs{i} + \abs{j} + \abs{k}}$~. Let

\[
U_V := \xi^i \partial_{x^i} - \xi^j \xi^i \rund{\Gamma_{i j}^{V, k} \circ \Pi_{\widehat{sT V^{|n}} }} \partial_{\xi^k} \in \fX\rund{\widehat{sT V^{|n}} }
\]

and

\[
U_W := \eta^i \partial_{y^i} - \eta^j \eta^i \rund{\Gamma_{i j}^{W, k} \circ \Pi_{\widehat{sT W^{|q}} }} \partial_{\eta^k} \in \fX\rund{\widehat{sT W^{|q}} } \,~.
\]

If $\Xi$ is affine w.r.t. $\nabla^V$ and $\nabla^W$ then $U_V$ and $U_W$ are $\widetilde \Xi$-related. Conversely, if $\nabla^V$ and $\nabla^W$ are torsionfree and $U_V$ and $U_W$ are $\widetilde \Xi$-related then $\Xi$ is affine w.r.t. $\nabla^V$ and $\nabla^W$~.
\end{lemma}

{\it Proof:} $\Xi$ being affine means $(d \Xi) \rund{\nabla^U_X Y} = \rund{\Xi^* \nabla^W}_X ((d \Xi) Y)$ for all $\P$- super vectorfields $X$ and $Y$ on $U^{|n}$~. The left hand side can be computed as

\[
(d \Xi) \rund{\nabla^U_X Y} = X^i \rund{\partial_i Y^k + (- 1)^{\abs{i} (\abs{Y} + \abs{j}) } Y^j \Gamma_{i j}^{V, k}} \rund{\partial_k \Xi^l} \rund{\Xi^* \partial_l} \,~,
% &=& (d \Xi) \rund{X^i \rund{\partial_i Y} + (- 1)^{\abs{i} (\abs{Y} + \abs{j}) } X^i Y^j \Gamma_{i j}^{V, k} \partial_k}
\]

and the right hand side with the help of the Christoffel symbols $\widehat \Gamma_{i k}^l = \rund{\partial_i \Xi^r} \rund{\Gamma_{r k}^{W, l} \circ \Xi}$ of $\Xi^* \nabla^W$ as

\begin{eqnarray*}
&& \rund{\Xi^* \nabla^W}_X ((d \Xi) Y) \\
% && \phantom{12} = X^i \rund{\Xi^* \nabla^W}_{\partial_i} \rund{Y^j \rund{\partial_j \Xi^k} \rund{\Xi^* \partial_k}} \\
&& \phantom{12} = X^i \rund{\rund{\partial_i Y^k} \rund{\partial_k \Xi^l} + (- 1)^{\abs{i} (\abs{Y} + \abs{j}) } Y^j \rund{\partial_i \partial_j \Xi^l} + (- 1)^{\abs{i} (\abs{Y} + \abs{k}) } Y^j \rund{\partial_j \Xi^k} \widehat \Gamma_{i k}^l} \rund{\Xi^* \partial_l} \\
&& \phantom{12} = X^i \left(\rund{\partial_i Y^k} \rund{\partial_k \Xi^l} + (- 1)^{\abs{i} (\abs{Y} + \abs{j}) } Y^j \rund{\partial_i \partial_j \Xi^l} \right. \\
&& \phantom{12345} \left. + (- 1)^{\abs{i} (\abs{Y} + \abs{k}) } Y^j \rund{\partial_j \Xi^k} \rund{\partial_i \Xi^r} \rund{\Gamma_{r k}^{W, l} \circ \Xi}\right) \rund{\Xi^* \partial_l} \,~.
\end{eqnarray*}

Therefore $\Xi$ is affine iff

\[
\Gamma_{i j}^{V, k} \rund{\partial_k \Xi^l} = \rund{\partial_i \partial_j \Xi^l} + (- 1)^{\abs{i} (\abs{j} + \abs{k}) } \rund{\partial_j \Xi^k} \rund{\partial_i \Xi^r} \rund{\Gamma_{r k}^{W, l} \circ \Xi} \,~.
\]

Now $U_V$ and $U_W$ $\widetilde \Xi$-related means $\rund{d \widetilde \Xi} U_V = {\widetilde \Xi}^* U_W$~. We drop the notations $\circ \Pi_{\widehat{sT V^{|n}} }$ and $\circ \Pi_{\widehat{sT W^{|q}} }$ for a better reading and compute the left hand side as

\[
\rund{d \widetilde \Xi} U_V = \xi^i \rund{\partial_i \Xi^l} \rund{{\widetilde \Xi}^* \partial_{y^l}} + \xi^j \xi^i \rund{\partial_i \partial_j \Xi^l - \Gamma_{i j}^{V, k} \rund{\partial_k \Xi^l}} \rund{{\widetilde \Xi}^* \partial_{\eta^l}} \,~,
\]

the right hand side as

\[
{\widetilde \Xi}^* U_W = \xi^i \rund{\partial_i \Xi^l} \rund{{\widetilde \Xi}^* \partial_{y^l}} - (- 1)^{\abs{i} (\abs{j} + \abs{k}) } \xi^j \xi^i \rund{\partial_j \Xi^k} \rund{\partial_i \Xi^r} \rund{\Gamma_{r k}^{W, l} \circ \Xi} \rund{{\widetilde \Xi}^* \partial_{\eta^l}} \,~.
\]

Therefore $U_V$ and $U_W$ are $\Xi$-related iff

\begin{eqnarray*}
&& 2 \partial_i \partial_j \Xi^l - \rund{\Gamma_{i j}^{V, k} + (- 1)^{\abs{i} \abs{j}} \Gamma_{j i}^{V, k}} \rund{\partial_k \Xi^l} \\
%&& \phantom{12} =  - (- 1)^{\abs{i} (\abs{j} + \abs{k}) } \rund{\partial_j \Xi^k} \rund{\partial_i \Xi^r} \rund{\Gamma_{r k}^{W, l} \circ \Xi} - (- 1)^{\abs{j} \abs{k}} \rund{\partial_i \Xi^k} \rund{\partial_j \Xi^r} \rund{\Gamma_{r k}^{W, l} \circ \Xi} \,~. \\
%&& \phantom{12} =  - (- 1)^{\abs{i} (\abs{j} + \abs{k}) } \rund{\partial_j \Xi^k} \rund{\partial_i \Xi^r} \rund{\Gamma_{r k}^{W, l} \circ \Xi} - (- 1)^{\abs{j} \abs{r}} \rund{\partial_i \Xi^r} \rund{\partial_j \Xi^k} \rund{\Gamma_{k r}^{W, l} \circ \Xi} \,~. \\
&& \phantom{12} =  - (- 1)^{\abs{i} (\abs{j} + \abs{k}) } \rund{\partial_j \Xi^k} \rund{\partial_i \Xi^r} \rund{\rund{\Gamma_{r k}^{W, l} + (- 1)^{\abs{r} \abs{k}} \Gamma_{k r}^{W, l}}  \circ \Xi} \,~,
\end{eqnarray*}

and the statement has become obvious. $\Box$ \\

Obviously the construction of $U$ and so also of $\Phi$ commute with ${}^\#$ and ${}^\rho$~.

\begin{cor} [Invariance of $\Phi$ under affine $\P$-supermorphisms] \label{geodcompatibleaffine} Let $\nabla^\N$ be a $\P$-connection on $\N$~.
\item[(i)] Let $\Xi: \M \rightarrow \N$ be a $\P$-supermorphism and $U_\N$ the generator of the geodesic flow on $\widehat{sT \N}$ w.r.t. $\nabla^\N$~. If $\Xi$ is affine w.r.t. $\nabla$ and $\nabla^\N$ then $U$ and $U_\N$ are $\widetilde \Xi$-related. If $\nabla$ and $\nabla^\N$ are torsionfree then also the converse is true.
\item[(ii)] There exists a unique torsionfree affine $\P$-connection $\bar \nabla$ on $\M$ such that $U$ is also the generator of the geodesic flow to $\bar \nabla$~. It is given by $\bar \nabla := \nabla - \frac{1}{2} T_\nabla$~.
\end{cor}

{\it Proof:} (i) obvious from lemma \ref{geodcompatibleaffinelocal}.

(ii) In local super charts define $\bar \nabla$ by the Christoffel symbols $\frac{1}{2} \rund{\Gamma_{i j}^k + (- 1)^{\abs{i} \abs{j}} \Gamma_{j i}^k}$~. Then obviously $U$ is the generator of its geodesic flow. Therefore also $\bar \nabla$ is globally defined by lemma \ref{geodcompatibleaffinelocal}. $\Box$ \\

Let us discuss two types of specializations of $\Phi$~. The first one are the integral curves to $U$ :

\begin{prop} There is a 1-1-correspondence between integral curves $\eta: I \rightarrow \widehat{sT \M}$ to $U$ and $\P$-curves $\gamma: I \rightarrow \M$ having $\nabla_{\dot \gamma} \dot \gamma = 0$ given by $\eta \mapsto \Pi_{\widehat{sT \M}} \circ \eta$ and $\dot \gamma \mapsfrom \gamma$~.
\end{prop}

{\it Proof:} In local super coordinates $\nabla_{\dot \gamma} \dot \gamma = 0$ is equivalent to the second order system of differential equations

\[
\ddot \gamma^k = - \dot \gamma^j \dot \gamma^i \rund{\Gamma_{i j}^k \circ \gamma} \,~,
\]

and $\eta = \rund{\sigma^i, \tau^k}$ being an integral curve to $U$ to the first order system of differential equations

\[
\dot \sigma^k = \tau^k \,~, \, \dot \tau^k = - \tau^j \tau^i \rund{\Gamma_{i j}^k \sigma} \,~.
\]

So the statement is obvious. $\Box$

\begin{defin} A $\P$-curve $\gamma: I \rightarrow \M$ is called a geodesic in $\M$ w.r.t. $\nabla$ iff $\nabla_{\dot \gamma} \dot \gamma = 0$~.
\end{defin}

\begin{cor} \label{largestgeod} Let $x \in \M^\P$ and $v \in \rund{sT_x \M}_0$~. Then there exists a largest geodesic \\
$\gamma: I \rightarrow \M$ such that $\gamma(0) = x$ and $\dot \gamma(0) = v$~. It is given by $\gamma = \Pi_{\widehat{sT \M}} \circ \Phi(\diamondsuit, X)$ with $I := I_{v^\#, U^\#}$~.
\end{cor}

Observe that for $\Q$ large enough $U$ and $U_\N$ being $\widetilde \Xi$-related is also equivalent to $\Xi$ mapping $\Q$-geodesics w.r.t. $\nabla$ to $\Q$-geodesics w.r.t. $\nabla^\N$~. \\

Here another nice characterization of geodesics:

\begin{theorem} A $\P$-curve $\gamma: I \rightarrow \M$ is a geodesic iff it is parallel w.r.t. the trivial affine connection $\nabla^\triv$ on $I$ and $\nabla$~.
\end{theorem}

{\it Proof:} The second means $(d \gamma) \rund{\nabla^{\triv}_X Y} = \rund{\gamma^* \nabla}_X ((d \gamma) Y)$ for all $X, Y \in \P \otimes \fX(I)$~. Let us write $X = f \partial_t$ and $Y = g \partial_t$~, $f, g \in \P \otimes \C^\infty(I)$~. Then the left hand side can be computed as $(d \gamma) \rund{\nabla^{\triv}_X Y} = f \dot g \dot \gamma$~, the right hand side as

\[
\rund{\gamma^* \nabla}_X \rund{(d \gamma) Y} = \rund{\gamma^* \nabla}_{f \partial_t} \rund{(d \gamma) g \partial_t} = f \dot g \dot \gamma + \nabla_{\dot \gamma} \dot \gamma \,~. \, \Box
\]

\begin{theorem} \label{velocity}
Denote by $\cdot: \rz \times \widehat{sT \M} \rightarrow \widehat{sT \M}$ the fiberwise multiplication in local super coordinates given by $\rund{x^i, t \xi^k}$~. Then

\[
\widetilde \Omega := \schweif{\left.(s, t, v) \in \rz^2 \times T M \, \right| \, (s, t \cdot v) \in \Omega} = \schweif{\left.(s, t, v) \in \rz^2 \times T M \, \right| \, (s t, v) \in \Omega} \,~,
\]

and

\[
\Phi \circ \rund{t, s \cdot \Pr\nolimits_{\widehat{sT \M}}} = s \cdot \Phi \circ \rund{s t, \Pr\nolimits_{\widehat{sT \M}}}
\]

as $\P$-supermorphisms from $\left.\rund{\rz \times \rz \times \widehat{sT \M}}\right|_{\widetilde \Omega}$ to $\widehat{sT \M}$~.
\end{theorem}

More concretely this says that if $\gamma$ is the largest geodesic with $\dot \gamma(0) = v \in \rund{sT_x \M}_0$ then $\eta := \gamma(s \diamondsuit)$ is the largest geodesic with $\dot \eta(0) = s \cdot v$~. \\

{\it Proof:} Let $s \in \rz$~. Then

\[
\varphi := \Phi \circ (\diamondsuit, s \cdot v): I_{s \cdot v^\#, U^\#} \rightarrow \widehat{sT \M}
\]

is an integral curve to $U$ with $\varphi(0) = s \cdot v$~. Also

\[
\eta := s \cdot \Phi \circ (s \diamondsuit, v): I \rightarrow \widehat{sT \M} \,~,
\]

$I := \schweif{t \in \rz \, \left| \, s t \in I_{v^\#, U^\#}\right.}$~, is an integral curve to $U$ with $\eta(0) = s \cdot v$~, indeed: In a local super chart write $\Phi(\diamondsuit, v) = (\sigma, \tau)$~, which is an integral curve to $U$~. So $\eta := (\sigma, s \tau)(s \diamondsuit)$~, and so

\[
\dot \eta = \rund{s \dot \sigma(\diamondsuit s), s^2 \dot \tau} = \rund{s \tau, - s^2 \tau^j \tau^i \rund{\Gamma_{i j}^k \circ \sigma}}(s \diamondsuit) = \eta^* U \,~.
\]

Therefore $I \subset I_{s \cdot v^\#, U^\#}$~. If $s = 0$ then $I = \rz$~. If $s \not= 0$ then by the same trick as used for $\eta$

\[
\vartheta := \frac{1}{s} \cdot \varphi\rund{\frac{1}{s} \diamondsuit}: s \ I_{s \cdot v^\#, U^\#} \rightarrow \widehat{sT \M}
\]

is an integral curve to $U$ with $\vartheta(0) = v$~. Therefore $s \ I_{s \cdot v^\#, U^\#} \subset I_{v^\#, U^\#}$~, and so $I = I_{s \cdot v\#, U^\#}$~. This gives the first statement. Furthermore $\varphi = \eta$~, and so we have also the second statement. $\Box$ \\

The second type of specializations which I would like to discuss are the geodesic super exponential maps:

\begin{defin}
Let $p \in \M^\P$~. The $\P$-supermorphism

\[
\exp_p := \Pi_{\widehat{sT \M}} \circ \Phi|_{\{1\} \times \widehat{sT_p \M}} : \left.\widehat{sT_p \M}\right|_{\Omega_{p^\#} } \rightarrow \M \,~,
\]

$\Omega_{p^\#} := \rund{\{1\} \times T_{p^\#} M} \cap \Omega$~, is called the geodesic exponential map to $p$ w.r.t. $\nabla$~.
\end{defin}

\begin{lemma}
\item[(i)] Let $x \in \M^\P$ and $v \in \rund{sT_x \M}_0$~. Then $\schweif{t \in \rz \, \left| \, t v^\# \in \Omega_{p^\#}\right.} = I_{v^\#, U^\#}$~, and \\
$\gamma := \exp_p(t v): I_{v^\#, U^\#} \rightarrow \M$ is the largest geodesic with $\dot \gamma(0) = v$~.
\item[(ii)] Identifying $sT_p \M \simeq sT_0 \ \widehat{sT_p \M}$ we have $\rund{d \exp_p}(0) = \id_{sT_p \M}$~, and so $\exp_p$ is a $\P$-superdiffeomorphism locally at $0$~.
\end{lemma}

{\it Proof:} (i) The first statement is an easy exercise using theorem \ref{velocity}. For the second statement by theorem \ref{velocity}

\[
\exp_p(t v) = \Pi_{\widehat{sT \M}} \circ \Phi(1, t v) = \Pi_{\widehat{sT \M}} \circ \rund{t \cdot \Phi(t, v)} = \Pi_{\widehat{sT \M}} \circ \Phi(t, v) \,~,
\]

which is the largest geodesic with $\dot \gamma(0) = v$ by corollary \ref{largestgeod}.

(ii) $\rund{d \exp_p}(0) v = v$ for all $v \in \rund{sT_p \M}_0$ by (i), but by passing from $\P$ to $\P \boxtimes \bigwedge \rz$ even for all $v \in sT_p \M$~. The rest follows by the super inverse function theorem. $\Box$ \\

As in \cite{Goertsches} proposition 4.11 for isometries between Riemannian supermanifolds and proposition 4.14 for Killing super vectorfields we obtain the faithful linearization result for affine $\P$-supermorphisms:

\begin{cor}[Faithful linearization] \label{linearize} Let $M$ be connected and $p \in \M^\P$~.
\item[(i)] Any affine $\P$-supermorphism $\Xi: \M \rightarrow \N$ is uniquely determined by $\Xi(p) \in \N$ and $(d \Xi)(p): sT_p \M \rightarrow sT_{\Xi(p)} \N$~.
\item[(ii)] Any $X \in \fX_\infaff(\M)$ is uniquely determined by $X(p)$ and $(\nabla X)(p)$~.
\item[(iii)] Let $\M$ and $\N$ be $\P$-Riemannian and of same super dimension. Then also any isometry $\Xi: \M \rightarrow \N$ is uniquely determined by $\Xi(p) \in \N$ and $(d \Xi)(p): sT_p \M \rightarrow sT_{\Xi(p)} \N$~.
\end{cor}

{\it Proof:} (i) Without restriction assume $\M$ is a supermanifold and $p \in M$~. Then we can use the canonical embedding $\rund{\id_M, {}^{\#}}: M \hookrightarrow \M$~. Assume also $\Sigma: \M \rightarrow \N$ is affine with $\Sigma(p) = \Xi(p)$ and $(d \Sigma)(p) = (d \Xi)(p)$~. Define

\[
U := \schweif{x \in M \, \left| \, \Sigma(x) = \Xi(x) \text{ and } (d \Sigma)(x) = (d \Xi)(x)\right.} \,~.
\]

Obviously $U$ is closed, and $p \in U$~. Now let $x \in U$ be arbitrary. Then there exists an open neighbourhood $V \subset M$ of $x$ such that $\Sigma|_{\M|_V} = \Xi_{\M|_V}$~. Indeed: Since $\Sigma$ and $\Xi$ are affine we obtain by corollary \ref{geodcompatibleaffine}

\[
\Sigma \circ \exp_x = \exp_{\Sigma(x)} \circ \widehat{(d \Sigma)(x)} = \exp_{\Xi(x)} \circ \widehat{(d \Xi)(x)} = \Xi \circ \exp_x \,~,
\]

but $\exp_x$ locally at $0$ is a $\P$-superdiffeomorphism.

Therefore $U$ is also open, so $U = M$ and $\Omega = \Id_\M$~.

(ii) Let $X \in \fX_\infaff\rund{\M}$ such that $X(p) = 0$ and $(\nabla X)(p) = 0$~. After passing from $\P$ to $\P \boxtimes \bigwedge \rz$ and from $X$ to $X_0 + \alpha X_1$~, $\alpha$ generating $\bigwedge \rz$~, we may assume $X$ even. For $\Q$ large enough take $a \in \Q_0 \setminus \{0\}$ such that $a^2 = 0$~. Then $J_X(a, \diamondsuit): \M \rightarrow \M$ is an affine $\Q$-superdifferomorphism by theorem \ref{Killingisom}, in local super charts given by $\Id_\M + a X$~. Therefore $J_X(a, p) = p$ and $\rund{d J_X(a, \diamondsuit)}(p) = \id_{sT_p \M}$~: indeed, since $X(p) = 0$ and $(\nabla X)(p) = 0$

\[
\rund{d J_X(a, \diamondsuit)}(p) \partial_i = \partial_i + a \rund{\rund{\partial_i X^k}(p)} \partial_k = \partial_i + \rund{\nabla_{\partial_i} X}(p) = \partial_i \,~.
\]

Therefore $J_X(a, \diamondsuit) = \Id_\M$ by (i), and so $X = 0$~.

(iii) direct consequence of (i) since by lemma \ref{isometric} (ii) $\Xi$ is affine w.r.t. the Levi-Civita connections. $\Box$ \\

{\bf From now on till the end of this article let $g$ be Riemannian $\P$-supermetric on $\M$~.} Then, as already pointed out in section \ref{vb}, it induces an isomorphism

\[
\varphi: sT \M \mathop{\rightarrow}\limits^{\sim} sT^* \M \,~, \, X \mapsto g(\phantom{1}, X) \,~.
\]

In the $\P$-supermanifold realization of $sT \M$ and $sT^* \M$~, $g$ becomes a function

\[
\widehat g \in \C^\infty\rund{\widehat{sT \M} \times_\M \widehat{sT \M}}_0 \,~,
\]

uniquely determined by the property that $\widehat g \circ \rund{\widehat X, \widehat Y} = g(X, Y)$ for all even $\P$- super vectorfields $X, Y$ on $\M$ and $\varphi$ a strong $\P$- super family isomorphism

\[
\begin{array}{ccc}
\widehat{sT \M} & \mathop{\rightarrow}\limits^{\widehat \varphi} & \widehat{sT^* \M} \\
{}_{\Pi_{\widehat{sT \M}}} \searrow & \circlearrowleft & \swarrow {}_{\Pi_{\widehat{sT^* \M}}} \\
& \M &
\end{array} \,~.
\]

In local super coordinates and using the standard frames $\rund{\partial_k}$ of $sT \M$ and $d x^k = \partial_k^*$ of $sT^* \M$~, $\widehat g$ is given by $\xi^i \eta^j \rund{g_{j i} \circ \Pi_{\widehat{sT \M} \times_\M \widehat{sT \M}} }$~, $g_{i j} := g\rund{\partial_i, \partial_j} \in \C^\infty(\M)_{\abs{i} + \abs{j}}$~, $\widehat \varphi$ by $\rund{\Id_\M, (- 1)^{\abs{i} + \abs{j}} \rund{g_{i j} \circ \Pi_{\widehat{sT \M}} } \xi^j}$ and ${\widehat \varphi}^{- 1}$ by $\rund{\Id_\M, (- 1)^{\abs{k} + \abs{l}} \rund{g^{k l} \circ \Pi_{\widehat{sT^* \M}} } a_k}$~, where $\rund{g^{k l}} \in \C^\infty(\M)_{\abs{i} + \abs{j}}$ are defined by $g_{i j} g^{k j} = \delta_i^k$ or equivalently $g^{i j} g_{i k} = \delta_k^j$~. Observe that $g^{l k} = (- 1)^{\abs{k} + \abs{l} + \abs{k} \abs{l}} g^{k l}$~.

\begin{theorem}[Super Gauß lemma]
Let $g$ be a Riemannian $\P$-supermetric on $\M$~, $p \in \M$ and $\exp_p: \left.\widehat{sT_p \M}\right|_{\Omega_p} \rightarrow \M$ the geodesic exponential map w.r.t. the Levi-Civita connection to $g$~. Then

\[
g\rund{\exp_p\rund{\widehat v}} \rund{\rund{d \exp_p}\rund{\widehat v} \diamondsuit, \rund{d \exp_p}\rund{\widehat v} v} = g(p)(\diamondsuit, v)
\]

on $sT_p \M$ for all $v \in \rund{sT_p \M}_0$ such that $v^\# \in \Omega_{p^\#}$~.
\end{theorem}

As in the classical case this implies that all geodesics emanating from $p \in M$ are perpendicular to all pseudo spheres $\exp_p \rund{\left.\schweif{\widehat{g(p)}(w, w) = c}\right|_{\Omega_{p^\#}} }$~, $c \in \P_0$~, centered at $p$~. \\

{\it Proof:} Let $w \in \rund{\Q \boxtimes_\P sT_p \M}_0$~. Then there exists $\eps > 0$ such that $s \rund{v^\# + t w^\#} \in \Omega_{p^\#}$ for all $s \in \ ]- \eps, 1 + \eps[$ and $t \in \ ]- \eps, \eps[$~. So we can define the $\Q$-supermorphism

\[
\varphi:= \exp_p(s (v + t w)): \ ]- \eps, 1 + \eps[ \ \times \ ]- \eps, \eps[ \ \rightarrow \ \M
\]

and the superfunctions

\[
f := \rund{\varphi^* g}\rund{(d \varphi) \partial_s, (d \varphi) \partial_s} \,~, \, h := \rund{\varphi^* g}\rund{(d \varphi) \partial_s, (d \varphi) \partial_t} \in \Q_0 \otimes \C^\infty\rund{]- \eps, 1 + \eps[ \ \times \ ]- \eps, \eps[} \,~.
\]

Since $\varphi(\diamondsuit, t)$ is a geodesic for all $t \in \ ]- \eps, \eps[$ and $\nabla$ and $g$ are compatible

\[
\partial_s f = 2 \rund{\varphi^* g}\rund{\rund{\varphi^* \nabla}_{\partial_s} (d \varphi) \partial_s, (d \varphi) \partial_s} = 0 \,~,
\]

and so $f = f(0, t) = g(p)(v + t w, v + t w)$~. Again since $\varphi(\diamondsuit, t)$ is a geodesic for all $t \in \ ]- \eps, \eps[$~, $\nabla$ and $g$ are compatible and $\nabla$ is torsionfree

\begin{eqnarray*}
\partial_s h &=& \rund{\varphi^* g}\rund{(d \varphi) \partial_s, \rund{\varphi^* \nabla}_{\partial_s} \rund{(d \varphi) \partial_t}} \\
&=& \rund{\varphi^* g}\rund{(d \varphi) \partial_s, \rund{\varphi^* \nabla}_{\partial_t} \rund{(d \varphi) \partial_s}} \\
&=& \frac{1}{2} \partial_t f = g(p)(v, w) \,~.
\end{eqnarray*}

Therefore

\[
g\rund{\exp_p\rund{\widehat v}} \rund{\rund{d \exp_p}\rund{\widehat v} \diamondsuit, \rund{d \exp_p}\rund{\widehat v} v} = h(1, 0) = h(0, 0) + g(p)(v, w) = g(p)(v, w) \,~. \, \Box
\]

\section{Mechanics on $\mathcal{P}$- Riemannian supermanifolds}

We have an even $\P$-linear map $\widetilde{\phantom{1}}: sT \M \rightarrow \rund{\Pi_{T M}}_* sT\rund{\widehat{sT \M}}$ assigning to every $\P$- super vectorfield $X$ on $\M$ the $\P$- super vectorfield $\widetilde X$ on $\widehat{sT \M}$ in local super coordinates of $\M$ given by $\widetilde X := \rund{X^i \circ \Pi_{\widehat{s T \M}} } \partial_{\xi^i}$ with $X = X^i \partial_i$~. Furthermore we have the even $\P$-linear sheaf morphism $\varphi^{- 1} \circ d: \C^\infty_\M \rightarrow sT \M$ assigning to every $f$ its super gradient vectorfield $X_{g, f}$ uniquely determined by the property $d f = g\rund{\diamondsuit, X_{g, f}}$~. \\

We define the kinetic energy function

\[
\T := \frac{1}{2} g \circ \Delta_\M \in \C^\infty\rund{\widehat{sT \M}}_0 \,~,
\]

where $\Delta_\M: \widehat{sT \M} \hookrightarrow \widehat{sT \M} \times_\M \widehat{sT \M}$ denotes the diagonal embedding, in local super coordinates given by $\rund{\Id_\M, \xi_k, \xi_k}$~. In local super coordinates $\T = \frac{1}{2} \xi^i \xi^j \rund{g_{j i} \circ \Pi_{\widehat{sT \M}} }$~, and $\T \circ \varphi^{-1}$ by $\frac{1}{2} (- 1)^{\abs{l}} a_k a_l \rund{g^{k l} \circ \Pi_{\widehat{sT^* \M}} }$~. Interpreting $\xi^i$ as velocity coordinates, $a_k$ become precisely the associated momentum coordinates since $a_k \circ \varphi = \partial_{\xi^k} \T$~. \\

Now assume $(\N, \omega)$ is a symplectic $\P$-supermanifold, so $\omega \in sT^* \N \boxtimes sT^* \N$ is a non-degenerate even $\P$-$2$-form, so graded antisymmetric. Then as for the super Riemannian case we have an even $\P$-linear sheaf morphism $\C^\infty_\N \rightarrow sT \N$ assigning to every $f$ its super Hamiltonian vectorfield $X_{\omega, f}$ uniquely determined by $\omega\rund{\phantom{1}, X_{\omega, f}} = d f$~. \\

Fortunately, $\widehat{sT^* \M}$ is symplectic in a canonical way: we have the canonical Liouville form $\alpha \in sT^*\rund{\widehat{sT^* \M}}$ in local super coordinates given by $\alpha = a_i d x^i$~. $\omega := d \alpha$ is non-degenerate since in local super coordinates $\omega = d a_i \wedge d x^i$~. In particular \\
$\omega\rund{\partial_{x^i} \otimes \partial_{a_k}} = - \delta_i^k = - (- 1)^{\abs{i}} \omega\rund{\partial_{a_k} \otimes \partial_{x^i}}$~. So for all $f \in \C^\infty\rund{\widehat{sT^* \M}}$

\[
X_{\omega, f} = \rund{\partial_{a_i} f} \partial_{x^i} - (- 1)^{\abs{r}} \rund{\partial_{x^r} f} \partial_{a_r} \,~.
\]

Finally let $\V \in \C^\infty(\M)_0$~. We regard $\V$ as a potential given on $\M$~. We take the Levi-Civita connection $\nabla$ on $\M$ w.r.t. $g$ according to theorem \ref{LeviCivita} and the generator of the geodesic flow on $U$ on $\widehat{sT \M}$ w.r.t. $\nabla$ from definition \ref{geodflow}. Then we have three descriptions of classical mechanics, which can be easily superized, each giving an answer to the question: when is a $\P$-curve $\gamma: I \rightarrow \M$ a trajectory of a `particle' moving in the potential $\V$ ?

\begin{itemize}
\item[(i)] Newton would say:~...iff $\nabla_{\dot \gamma} \dot \gamma = - \gamma^* X_{g, \V}$~, or equivalently iff $\dot \gamma$ is an integral curve to the $\P$- super vectorfield $U - \widetilde{X_{g, \V}}$ on $\widehat{sT \M}$~.

\item[(ii)] Lagrange would say:~...iff in local super coordinates

\begin{equation} \label{Lagrange}
\rund{\rund{\partial_{\xi^k} \L} \circ \dot \gamma}^\cdot = \rund{\partial_{x^k} \L} \circ \dot \gamma \,~,
\end{equation}

where $\L := \T - \V \circ \Pi_{\widehat{sT \M}} \in \C^\infty\rund{\widehat{sT \M}}_0$ denotes the Lagrangian, or equivalently $\int_K \rund{\L \circ \dot \gamma}(t) d t$ is stationary under local variation of the trajectory $\gamma|_K$ for all $K \subset I$ compact~.

\item[(iii)] Hamilton would say:~...iff $\varphi \circ \gamma$ is an integral curve to the Hamiltonian super vectorfield $X_{\omega, \H}$ on $\widehat{sT^* \M}$ of the Hamilton function

\[
\H := \T \circ \varphi^{- 1} + \V \circ \Pi_{\widehat{sT^* \M}} \in \C^\infty(sT^* \M)_0 \,~.
\]

\end{itemize}

All three descriptions are infact equivalent, as the last result of this article tells us:

\begin{theorem}
\item[(i)] $\dot \gamma: I \rightarrow \widehat{sT \M}$ is an integral curve to $U - \widetilde{X_{g, \V}}$ iff it fulfills (\ref{Lagrange}).
\item[(ii)] $U - \widetilde{X_{g, \V}}$ and $X_{\omega, \H}$ are related under $\widehat \varphi$~.
%, and so $\varphi$ transforms the geodesic flow $\Phi$ on $sT \M$ into the integral flow $J_{X_{\omega, \H}}$ to the Hamiltonian vectorfield $X_{\omega, \H}$ :

%\[
%\begin{array}{ccc}
%(\rz \times sT \M)|_\Omega & \mathop{\mathop{\longrightarrow}\limits^\sim}\limits^{(\id, \varphi)} & \left.\rund{\rz \times sT^* \M}\right|_\Omega \\
%\Phi \downarrow \phantom{\Phi} & \circlearrowleft & \phantom{J_{X_{\omega, \H}}} \downarrow J_{X_{\omega, \H}} \\
%sT \M & \mathop{\mathop{\longrightarrow}\limits_\sim}\limits_{\varphi} & sT^* \M
%\end{array} \,~.
%\]

\end{theorem}

{\it Proof:}
(i) The left hand side of (\ref{Lagrange}) can be computed as

\[
\rund{\rund{\partial_{\xi^r} \L} \circ \dot \gamma}^\cdot = \rund{\dot{\gamma^k} \rund{g_{k r} \circ \gamma}}^\cdot = \ddot{\gamma^k} \rund{g_{k r} \circ \gamma} + \dot{\gamma^j} \dot{\gamma^i} \rund{\rund{\partial_i g_{j r}} \circ \gamma} \,~,
\]

and with the help of the formulas (\ref{LeviCivitaform}) of theorem \ref{LeviCivita} and $\partial_r \V = X_\V^r g_{k r}$ the right hand side as

\begin{eqnarray*}
\rund{\partial_{x^r} \L} \circ \dot \gamma &=& \frac{(- 1)^{\abs{r} (\abs{j} + \abs{i}) }}{2} \dot{\gamma^j} \dot{\gamma^i} \rund{\rund{\partial_r g_{i j}} \circ \gamma} - \rund{\partial_r \V} \circ \gamma \\
&=& \dot{\gamma^j} \dot{\gamma^i} \rund{\rund{\frac{1}{2}\rund{\partial_i g_{j r} + (- 1)^{\abs{i} \abs{j}} \partial_j g_{i r}} - \Gamma_{i j}^k g_{k r}} \circ \gamma} - \rund{X_\V^r g_{k r}} \circ \gamma \\
&=& \dot{\gamma^j} \dot{\gamma^i} \rund{\rund{\partial_i g_{j r} - \Gamma_{i j}^k g_{k r}} \circ \gamma} - \rund{X_\V^r g_{k r}} \circ \gamma \,~.
\end{eqnarray*}

Therefore and since $g$ is non-degenerate (\ref{Lagrange}) is equivalent to

\[
\ddot{\gamma^k} = - \dot{\gamma^j} \dot{\gamma^i} \rund{\Gamma_{i j}^k \circ \gamma} - X_\V^r \circ \gamma \,~,
\]

which says that $\dot \gamma$ is an integral curve to $U - \widetilde{X_{g, \V}}$~.

(ii) We use local super coordinates of $\M$ and again drop the notations $\circ \Pi_{\widehat{sT \M}}$ and $\circ \Pi_{\widehat{sT^* \M}}$ for a better reading. Using the formulas

\[
\rund{d \widehat \varphi} \partial_{x^i} = {\widehat \varphi}^* \partial_{x^i} + (- 1)^{\abs{r} + \abs{j}} \rund{\partial_i g_{r j}} \xi^j \rund{{\widehat \varphi}^* \partial_{a_r}} \phantom{1}~, \phantom{1} \rund{d \widehat \varphi} \partial_{\xi^k} = (- 1)^{\abs{r} (\abs{k} + 1)} g_{r k} \rund{{\widehat \varphi}^* \partial_{a_r}} \phantom{1}~,
\]

(\ref{LeviCivitaform}) of theorem \ref{LeviCivita}, $\partial_r \V = X_\V^r g_{k r}$~,

\[
\partial_r g_{i j} = - (- 1)^{\abs{r} \rund{\abs{i} + \abs{l}} } g_{i l} \rund{\partial_r g^{k l}} g_{k m}
\]

and

\[
X_{\omega, \H} = (- 1)^{\abs{j}} a_j g^{i j} \partial_{x^i} - (- 1)^{\abs{r}} \rund{\frac{(- 1)^{\abs{r} (\abs{k} + \abs{l}) + \abs{l}} }{2} a_k a_l \rund{\partial_r g^{k l}} + \partial_r \V} \partial_{a_r}
\]

we obtain

\begin{eqnarray*}
\rund{d \widehat \varphi} \rund{U - \widetilde{X_{g, \V}}} &=& \xi^i \rund{d \widehat \varphi} \partial_{x^i} - \rund{\xi^j \xi^i \Gamma_{i j}^k + X_\V^k} \rund{d \widehat \varphi} \partial_{\xi^k} \\
&=& \xi^i \rund{{\widehat \varphi}^* \partial_{x^i} + (- 1)^{\abs{r} + \abs{j}} \rund{\partial_i g_{r j}} \xi^j \rund{{\widehat \varphi}^* \partial_{a_r}} } \\
&& \phantom{12} - \rund{\xi^j \xi^i \Gamma_{i j}^k + X_\V^k} \rund{(- 1)^{(\abs{k} + 1) \abs{r}} g_{r k} \rund{{\widehat \varphi}^* \partial_{a_r}} } \\
&=& \xi^i \rund{{\widehat \varphi}^* \partial_{x^i}} + (- 1)^{\abs{r}} \rund{\xi^j \xi^i \rund{\partial_i g_{j r} - \Gamma_{i j}^k g_{k r}} - X_\V^r g_{k r}} \rund{{\widehat \varphi}^* \partial_{a_r}} \\
&=& \xi^i \rund{{\widehat \varphi}^* \partial_{x^i}} + (- 1)^{\abs{r}} \rund{\frac{(- 1)^{\abs{r} (\abs{i} + \abs{j}) }}{2} \xi^j \xi^i \rund{\partial_r g_{i j}}  - \partial_r \V} \rund{{\widehat \varphi}^* \partial_{a_r}} \\
&=& \xi^q g_{q j} g^{i j} \rund{{\widehat \varphi}^* \partial_{x^i}} \\
&& \phantom{12} - (- 1)^{\abs{r}} \rund{\frac{(- 1)^{\abs{r} (\abs{j} + \abs{l}) }}{2} \xi^j \xi^i g_{i l} \rund{\partial_r g^{k l}} g_{k j} + \partial_r \V} \rund{{\widehat \varphi}^* \partial_{a_r}} \\
&=& (- 1)^{\abs{q}} g_{j q} \xi^q g^{i j} \rund{{\widehat \varphi}^* \partial_{x^i}}  \\
&& \phantom{12} - (- 1)^{\abs{r}} \rund{\frac{(- 1)^{\abs{r} (\abs{k} + \abs{l}) + \abs{i} + \abs{j} + \abs{k}} }{2} g_{k j} \xi^j g_{l i} \xi^i \rund{\partial_r g^{k l}} + \partial_r \V} \rund{{\widehat \varphi}^* \partial_{a_r}} \\
&=& {\widehat \varphi}^* \rund{(- 1)^{\abs{j}} a_j g^{i j} \partial_{x^i} - (- 1)^{\abs{r}} \rund{\frac{(- 1)^{\abs{r} (\abs{k} + \abs{l}) + \abs{l}} }{2} a_k a_l \rund{\partial_r g^{k l}} + \partial_r \V} \partial_{a_r}} \\
&=& {\widehat \varphi}^* X_{\omega, \H} \,~. \, \Box
\end{eqnarray*}

\end{document}